\documentclass{article}

\usepackage{arxiv}
\usepackage{graphicx}
\usepackage{amsthm,amsmath}
\usepackage[numbers]{natbib}
\usepackage[colorlinks,citecolor=blue,urlcolor=blue]{hyperref}

\usepackage[english]{babel} 
\usepackage[utf8]{inputenc}
\usepackage[T1]{fontenc}
\usepackage{textcomp}
\usepackage{dsfont}
\usepackage{mathtools}
\usepackage{eqnarray}
\usepackage{framed}
\usepackage{tabularx}
\usepackage{booktabs} 
\usepackage[para]{threeparttable} 
\usepackage{pifont}
\usepackage{tkz-euclide}
\usepackage{diagbox}
\usepackage{xcolor}
\usepackage{amssymb, amsfonts, bm}
\usepackage{float}
\usepackage{subcaption}
\usepackage{xfrac}
\usepackage{mathrsfs}

\DeclareMathOperator{\diag}{diag}

\newcommand{\cF}{\mathcal{F}_\mu}
\newcommand{\cR}{{\mathcal{R}_\mu}}
\newcommand{\cE}{{\mathcal{E}}} 

\newcommand{\cY}{\mathcal{Y}}

\newcommand{\cG}{\mathcal{G}}

\renewcommand{\S}{\mathcal{S}^2}

\newcommand{\Sd}{{\mathcal{S}^{d-1}}}
\newcommand{\Sdmtwo}{\mathcal{S}^{d-2}}

\newcommand{\ts}{T_\mu\mathcal{S}^2}
\newcommand{\tsd}{T_\mu\mathcal{S}^{d-1}}

\newcommand{\Rd}{\mathbb{R}^d}

\newcommand{\Rdd}{\mathbb{R}^{d\times d}}

\newcommand{\R}{\mathbb{R}}

\newcommand{\one}{1}
\DeclareMathOperator{\minor}{minor}
\DeclareMathOperator{\major}{major}

\numberwithin{equation}{section}
\theoremstyle{plain}

\newtheorem{definition}{Definition}[section]

\newtheorem{proposition}{Proposition}[section]

\title{Quantiles and depth for directional data from elliptically symmetric distributions}

\author{ Konstantin Hauch \\
	Department of Mathematics,\\ 
    Technische Universität Kaiserslautern,\\ 
    Kaiserslautern, Germany\\
    \texttt{hauch@mathematik.uni-kl.de} \\
	\And
	Claudia Redenbach \\
	Department of Mathematics,\\ 
    Technische Universität Kaiserslautern,\\ 
    Kaiserslautern, Germany\\
    \texttt{redenbach@mathematik.uni-kl.de} \\
}

\begin{document}

\maketitle 

\begin{abstract}
We present canonical quantiles and depths for directional data following a distribution which is elliptically symmetric about a direction $\mu$ on the sphere $\Sd$.
Our approach extends the concept of Ley et al. \cite{Ley14}, which provides valuable geometric properties of the depth contours (such as convexity and rotational equivariance) and a Bahadur-type representation of the quantiles.
Their concept is canonical for rotationally symmetric depth contours. 
However, it also produces rotationally symmetric depth contours when the underlying distribution is not rotationally symmetric. 
We solve this lack of flexibility for distributions with elliptical depth contours.
The basic idea is to deform the elliptic contours by a diffeomorphic mapping to rotationally symmetric contours, thus reverting to the canonical case in Ley et al. \cite{Ley14}. 
A Monte Carlo simulation study confirms our results. 
We use our method to evaluate the ellipticity of depth contours and for trimming of directional data.
The analysis of fibre directions in fibre-reinforced concrete underlines the practical relevance.
\end{abstract}



\section{Introduction}

The classes of rotationally symmetric distributions and elliptically symmetric distributions in $\Rd$ have been well investigated by Kelker \cite{Kelker70}, Cambanis et al. \cite{Cambanis81} and Fang et al. \cite{Fang90}.
A random vector $V\in\Rd$ has a rotationally symmetric distribution if $V \overset{D}{=} OV$ for all $O\in SO(d)$ where $\overset{D}{=}$ refers to equality in distribution. 
Furthermore, every random vector $V\in\Rd$ following a rotationally symmetric distribution can be represented as $V\overset{D}{=}RU$ , where $U\sim Unif(\Sd)$ is independent of the real-valued random variable $R\sim F_R$. 
$U$ gives the direction of $V$ and $R$ is the length of $V$.  
Rotationally symmetric distributions are often regarded as the most natural non-uniform distributions in $\Rd$.  
For instance, the charge distribution of an electric field is rotationally symmetric around its source.
However, not all phenomena observed in practice can be represented by symmetric models.

Elliptically symmetric distributions extend the class of rotationally symmetric distributions.
The distribution of a random vector $W$ is elliptically symmetric if and only if $W\overset{D}{=}R\Sigma^{1/2}U$ with $U\sim Unif(\Sd)$, real-valued $R\sim F_R$ independent of $U$, and $\Sigma\in\Rdd$ a symmetric, positive definite matrix.
A random vector $W$ with an elliptically symmetric distribution can be transformed into a random vector $V=RU$ with a rotationally symmetric distribution via 
\begin{align}
    \Sigma^{-1/2}W &\overset{D}{=} R\Sigma^{-1/2} \Sigma^{1/2}U = R U=V.
    \label{eq:Elliptical symmetry is connected to rotational symmetry}
\end{align}

These concepts of symmetry transfer to the unit sphere $\Sd$, i.e., the case of directional data.
Distributions on $\Sd$ which are rotationally symmetric about a direction $\mu\in\Sd$ are also often regarded as the natural non-uniform distributions on $\Sd$ \cite{Verdebout20}. In most cases, rotationally symmetric distributions have tractable normalising constants.
Note that the density of a rotationally symmetric distribution is proportional to a function $f(x^T\mu)$. Thus, a projection onto $\mu$ enables a one-dimensional analysis of the distribution, for example, its concentration around $\mu$.
The class of distributions with rotational symmetry about $\mu\in\Sd$ is denoted by $\cR$.

In practice, symmetric models are often too restrictive. 
For instance, Leong and Carlile \cite{Leong1998} illustrated that rotational symmetry about a direction is a too strong assumption for a directional data set from neurosciences.
Kent \cite{Kent82} has fitted his elliptical model to a data set of 34 measurements of the directions of magnetisation for samples from the Great Whin Sill (Northumberland, England).
As in $\Rd$, distributions that are elliptically symmetric about a direction $\mu$ on $\Sd$ are an extension of the rotationally symmetric distributions.
The contours are more flexible to form elliptical shapes.
Due to the curved shape of the sphere, the transition from distributions which are rotationally symmetric about $\mu\in\Sd$ to distributions which are elliptically symmetric about $\mu$ is not obvious.
A remedy to this problem is to linearise $\Sd$ at some base point $\mu$ by considering the tangent space $\tsd$ at $\mu$. By using the theory for $\Rd$, a transformation between the two distributions can then be defined in the tangent space.

Ley et al. \cite{Ley14} introduced a concept of quantiles and depth for directional data.
They showed that their quantiles are asymptotically normal and established a Bahadur-type representation \cite{Bahadur66} for directional data $X\sim F_X \in \cR$. 
A Monte Carlo simulation study corroborated their theoretical results.
Statistical tools, like directional DD- and QQ-plots and a quantile-based goodness-of-fit test, were defined and illustrated on a data set of cosmic rays.
Their results are canonical for rotationally symmetric distributions. 
But their concept suffers from the disadvantage of producing rotationally symmetric depth contours, even if the underlying distribution has elliptical contours \cite{Pandolfo17}.

In this paper, we present a procedure solving the latter issue if the underlying distribution has elliptical contours.
The paper is organised as follows.
In Section \ref{sec:Basics-Quantile}, we first introduce basics about the distributions under consideration, extend the Mahalanobis transformation to $\Sd$, and summarise the findings of Ley at al. \cite{Ley14}.
Section \ref{sec:Quantiles-ellipse} contains our main contribution.
The idea is to map the unit vectors into the tangent space $\tsd$ where $\mu$ is the median direction of the observed sample. 
The mapped vectors are elliptically symmetric around the origin in $\tsd$. 
The multivariate Mahalanobis transformation \cite{Pennec06,Hardle03} is then applied in $\tsd$ to obtain a rotationally symmetric sample in $\tsd$. Mapping it back to $\Sd$,
we obtain a sample of unit vectors which are rotationally symmetric about $\mu$.
Thus, we can exploit the results from \cite{Ley14}.
All transformations are diffeomorphic such that we can trace back the results to the original unit vectors.
Section \ref{seq:Applications Quantiles} affirms our findings by a Monte Carlo study.
Furthermore, we apply our approach to a real-world data set from \cite{Maryamh20}: Directions of short steel fibres crossing a crack in ultra-high performance fibre-reinforced concrete (UHPFRC).

\section{Basics}\label{sec:Basics-Quantile}
\subsection{Rotational and elliptical symmetry about a direction on \texorpdfstring{$\Sd$}{Sd}}\label{sec:Rotational and elliptical symmetry about a direction}
\begin{definition}[Rotational symmetry about a direction]
	Let $X\in\Sd$ be a random vector and $\mu\in\Sd$. 
	The distribution of $X$ is rotationally symmetric about $\mu$ on $\Sd$ if and only if $X\overset{D}{=}OX$ for every $O\in SO(d)$ fulfilling $O\mu=\mu$. 
\end{definition}
Let $\cR$ be the class of distributions which are rotationally symmetric about $\mu\in\Sd$.
Projecting $X$ onto a vector space orthogonal to $\mu$ yields rotationally symmetric contours.
Distributions $F_X\in\cR$ are characterised by densities of the form
\begin{align}
	 f_{\mu}(x) = c_{d}f(x^T\mu), \quad x\in\Sd,
	\label{eq:cR-density}
\end{align}
where $f:[-1,1] \rightarrow\R_{\geq 0}$ is absolutely continuous and $c_{d}$ a normalising constant \cite{Ley14}.
The distribution of $X^T\mu$ is absolutely continuous w.r.t the Lebesgue measure on $[-1,1]$ \cite{Mardia99}.
The density of $X^T\mu$ reads
\begin{align}
	\Tilde{f}(t)=\omega_{d-1}c_{d}(1-t^2)^{\frac{d-3}{2}}f(t),
	\label{eq:Ytmu-density}
\end{align}
where $\omega_{d-1}$ is the surface area of $\Sdmtwo$ \cite{Verdebout20}.
A widely known distribution in $\cR$ is the von Mises-Fisher distribution where  $f(t)=\exp(\kappa t).$

\begin{definition}[von Mises-Fisher distribution $M_d(\mu,\kappa)$\cite{Mardia99}]\label{def:von Mises-Fisher distribution} 
	The probability density function of the von Mises–Fisher distribution is given by 
	\begin{equation}
		f_{vMF,\mu, \kappa}(x)=c_{d}\exp{(\kappa x^T\mu)},
		\label{eq:densityvMF}
	\end{equation}
	where $\kappa\geq 0$ is a concentration parameter, $\mu\in\Sd$ the mean direction, and $c_{d}$ the normalising constant.
\end{definition}
The concentration around $\mu$ increases with $\kappa$. 
The von Mises-Fisher distribution is unimodal for $\kappa >0$. For $\kappa=0$, we get the uniform distribution on the sphere. 

A generalisation of the von Mises-Fisher distribution is the Fisher-Bingham distribution \cite{Mardia99}, where a general quadratic equation is added in the exponent of the density in \eqref{eq:cR-density}. An example is the Kent distribution \cite{Kent82}.

\begin{definition}[Kent distribution $K(\mu,\kappa, A)$]\label{def:Kent distribution} 
	The probability density function of the Kent distribution is given by 
	\begin{equation}
		f_{K, \mu, A}(x)=c_{d}\exp{(\kappa x^T\mu + x^TAx)},
		\label{eq:densityKent}
	\end{equation}
	where $\kappa\geq 0$ is a concetration parameter, $\mu\in\Sd$ the mean direction, $A\in Sym(d)$ with $A\mu=0_d$ a shape parameter, and $c_{d}$ the normalising constant.
	The concentration around $\mu$ increases with $\kappa$, while $A\in Sym(d)$ controls the shape of the density contours.	 
\end{definition}
For large $\kappa$, the Kent distribution has a mode at $\mu$ and density contours which are elliptical \cite[p.177]{FLE87}.

Let $\cF$ be the class of distributions on $\Sd$ with a bounded density that admit a unique modal direction $\mu$.
We further assume that $\mu$ coincides with the Fisher spherical median \cite{Fisher85}, that is
\begin{align}
	\mu= \arg\min_{\gamma\in\Sd}E(\arccos(X^T\gamma)).   
	\label{eq:Fisher spherical median} 
\end{align}
For i.i.d. random vectors $X,X_1,\dots,X_n\in \Sd$ with $X\sim F\in\cF$, we estimate  $\mu$ by the root-$n$ consistent empirical Fisher spherical median \cite{Fisher85}
\begin{align}
	\hat{\mu}= \arg\min_{\gamma\in\S}\sum_{i=1}^N\arccos(X_i^T\gamma).     
	\label{eq:empirical Fisher spherical median}
\end{align}

Note that the definition of the class $\cR$ does not include that $\mu$ is the unique modal direction.
Here, we restrict attention to distributions in $\cR\cap\cF$ from now on. 
E.g., $M_d(\mu,\kappa)\in\cR\cap\cF$ for $\kappa>0$, and $K(\mu,\kappa, A)\in\cF$ for $\kappa>0$ under suitable conditions on $A\in Sym(d)$ given in the next section.

\subsection{Differential geometry}

Differential geometry examines smooth manifolds using differential and integral calculus as well as linear and multi-linear algebra.
It originates in studying spherical geometries related to astronomy and the geodesy of the earth.
For an introduction to differential geometry, see e.g. \cite{Tu11}.

We saw in \eqref{eq:Elliptical symmetry is connected to rotational symmetry} that a linear transformation $\Sigma$ transforms a random vector with a rotationally symmetric distribution into a random vector with an elliptically symmetric distribution in $\Rd$.
We want to proceed analogously for distributions on the sphere.
However, in general, the linear transformation $\Sigma$ does not necessarily map $\Sd$ onto itself. A remedy is provided by linearising the sphere $\Sd$ at a base point $\mu\in\Sd$.

The tangent space $\tsd$ to $\Sd$ at base point $\mu\in\Sd$ is the collection of all tangent vectors to $\Sd$ at $\mu$. 
It is a local Euclidean vector space with local origin in $\mu$.
Given $\mu\in\Sd$ and a tangent vector $v\in\tsd$, there is a unique geodesic from $\mu\in\Sd$ to some $x\in\Sd$ given as a mapping
\begin{align}
	c^{\mu,v}:[0,1]\rightarrow\Sd, 
	\label{eq:geodesic}   
\end{align}
starting at $c^{\mu,v}(0)=\mu$ with initial velocity $\dot c^{\mu,v}(0)=v$ and ending in $c^{\mu,v}(1)=x$ \cite{Fletcher08}.

In the following, we define mappings between the tangent space and the sphere.
\begin{definition}[Riemannian exponential map]
The Riemannian exponential map
\begin{align}
	Exp_\mu: \tsd \rightarrow \Sd   
	\label{eq:Exp_mu} 
\end{align}
maps a vector $v\in\tsd$ to $\Sd$ along the geodesic $c^{\mu,v}$ such that $x=Exp_\mu(v) =c^{\mu,v}(1)$. 
\end{definition}
The exponential map is locally diffeomorphic onto $V(\mu)= \Sd \setminus \{-\mu\}$, where $-\mu$ is called cut point and the set $\{-\mu\}$ is called cut locus.
Within $V(\mu)$ the exponential map $Exp_\mu$ has an inverse, the Riemannian logarithmic map.

\begin{definition}[Riemannian logarithmic map]
The Riemannian logarithmic map
\begin{align}
	Log_\mu: V(\mu) \rightarrow \tsd   
	\label{eq:Log_mu} 
\end{align}
maps a vector $x\in\Sd$ into $\tsd$ with $Exp_\mu(Log_\mu(x))=x$.
\end{definition}

The distance between $\mu$ and a point on the sphere is described by the Riemannian distance function.
\begin{definition}[Riemannian distance function]
For any point $x\in V(\mu)$, the Riemannian distance function is given by
\begin{align}
	d_{GD}(\mu, x) = ||Log_\mu(x)||_2=\arccos{(x^T\mu)} \in [0,\pi).    
	\label{eq:d_GD}
\end{align}
\end{definition}
Consequently, $Log_\mu(x)\in B_{d-1}(\pi) \subset\tsd$ with $B_{d}(r)$ the $d$-dimensional open ball of radius $r>0$ centred at the origin $0_d$.

\subsubsection{Examples on \texorpdfstring{$\S$}{S2}}
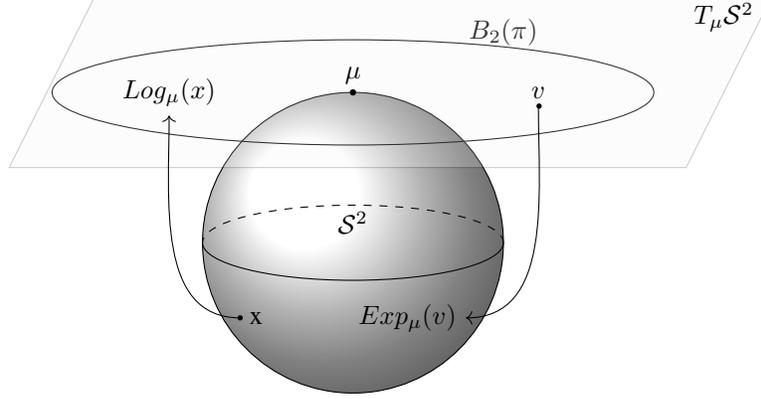
\begin{figure}
	\centering
	\begin{tikzpicture}[
        point/.style = {draw, circle, fill=black, inner sep=0.7pt},
        ]
        \def\rad{2cm}
        \coordinate (O) at (0,0); 
        \coordinate (N) at (0,\rad); 
        
        \filldraw[ball color=white] (O) circle [radius=\rad];
        \draw[dashed] 
          (\rad,0) arc [start angle=0,end angle=180,x radius=\rad,y radius=5mm];
        \draw
          (\rad,0) arc [start angle=0,end angle=-180,x radius=\rad,y radius=5mm];
        \draw
          (2*\rad,\rad) arc [start angle=0,end angle=360,x radius=2*\rad,y radius=7mm];
        \node at (\rad,2.8) {$B_{2}(\pi)$};  
        
        \begin{scope}[xslant=0.5,yshift=\rad,xshift=-2]
        \filldraw[fill=gray!10,opacity=0.2]
            (-5,1.25) -- (4,1.25) -- (4,-1) -- (-5,-1) -- cycle;
        \node at (3.5,1) {$T_\mu\mathcal{S}^2$};  
        \end{scope}
        
        \draw[dashed]
          (N) node[above] {$\mu$};
        \node[point] at (N) {};
        \draw[dashed]
          (O) node[above] {$\mathcal{S}^2$};
        
        \node[anchor=east] at (-1.7,\rad) (text) {$Log_\mu(x)$};
        \node[anchor=west] at (-1.5,-\rad/2) (description) {x};
        \draw (description) edge[out=180,in=-90,->] (text); 
        \node[anchor=east] at (1.5,-\rad/2) (text) {$Exp_\mu(v)$};
        \node[anchor=west] at (2.25,\rad) (description) {$v$};
        \draw (description) edge[out=-90,in=0,->] (text); 
        \fill (-1.5,-\rad/2) circle (1pt);
        \fill (2.474,\rad/1.1) circle (1pt);
	\end{tikzpicture}
	\caption{
		The tangent space $\ts$ of the sphere $\S$ and related operators.
	}
	\label{fig:sphere-tangent-space}
\end{figure}
The locally diffeomorphic exponential map on $\S$ reads 
\begin{align}
	Exp_\mu: \ts &\rightarrow \S\setminus\{-\mu\}, \nonumber\\
	v &\mapsto \mu \cos{(||v||_2)} + \frac{v}{||v||_2}\sin{(||v||_2)},
	\label{eq:defExp}
\end{align}
while the logarithmic map is given by
\begin{align}
	Log_\mu: \S \setminus\{-\mu\} &\rightarrow \ts, \nonumber \\
	x &\mapsto \frac{d_{GD}(x,\mu)}{\sin{(d_{GD}(x,\mu))}} z.
	\label{eq:defLog}
\end{align}
Here,  $z=(I_d-\mu\mu^T)x$ is the tangential part of $x$ and we use the convention $\frac{0}{\sin{(0)}}=1$ \cite{Hauberg18,Angulo12}.
See Figure \ref{fig:sphere-tangent-space} for an illustration.

With $\mu=(0,0,1)^T$ the logarithmic map in spherical coordinates $\phi \in [0, 2 \pi)$ and $\theta \in [0,\pi] $ reads 
\begin{align}
	Log_\mu(x)
	&
	=
	\left(
	\begin{array}{c}
		\cos{(\phi)}\sin{(\theta)}\\
		\sin{(\phi)}\sin{(\theta)}\\
		0\\
	\end{array}
	\right) \frac{\theta}{\sin{(\theta)}}
	= 
	\theta
	\left(
	\begin{array}{c}
		\cos{(\phi)}\\
		\sin{(\phi)}\\
		0\\
	\end{array}
	\right)
	\label{eq:Log_x_circle_ellipse}
	\in \ts,            
\end{align}
with $d_{GD}(x,\mu)=\arccos{(\cos(\theta))}=\theta$.
Thus, the length of any vector $Log_\mu(x)\in\ts$ coincides with the co-latitude angle $\theta$.
For a random vector $X$ with a distribution from $\cR \cap \cF$, we have $\Phi\sim Unif[0,2\pi]$ and $\Phi$ and $\Theta$  are independent. Hence,  
$Log_\mu(X)$ has circular density contours in $\ts$.

\subsection{The Mahalanobis transformation on \texorpdfstring{$\Sd$}{Sd}}
The idea of the Mahalanobis transformation in $\Rd$ is to linearly transform a real-valued data matrix into a centred, standardised and uncorrelated data matrix, see e.g. \cite{Hardle03}.
\begin{figure}
	\centering
	\includegraphics[width=0.8\linewidth]{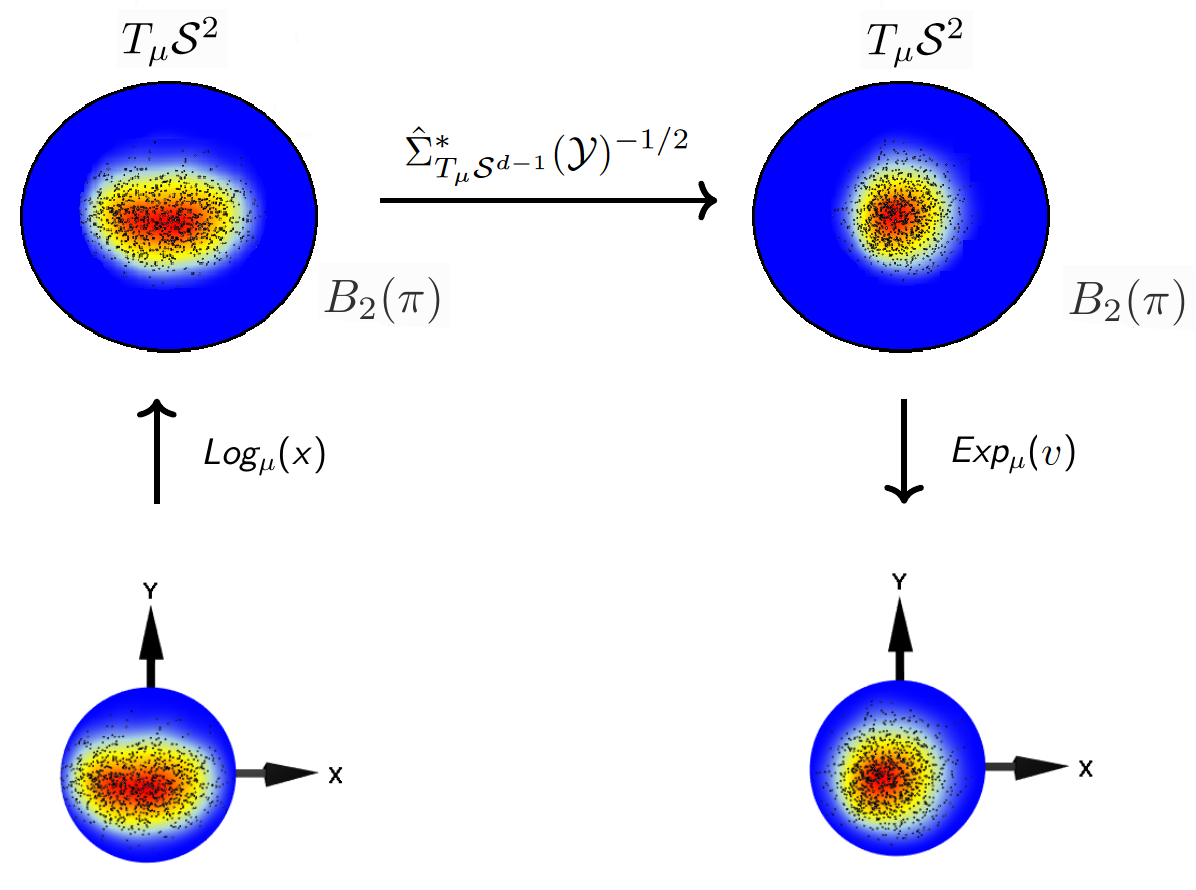}
	\caption
	{
		The Mahalanobis transformation of  $y_1,\dots,y_{500}\in\S$. The realisations are i.i.d. Kent distributed with $\mu=(0,0,1)$, $\kappa=12$ and $A = \diag(\beta,-\beta,0)$ with $\beta=5$.
		The Z-direction points out of the page.
	}
	\label{fig:mahalanobis-trafo-procedure}
\end{figure}

The Mahalanobis transformation can be generalised to the Riemannian manifold $\Sd$, see \cite{Pennec06}.
Here, the row vectors $y_1,\dots,y_n\in\Sd\setminus\{-\mu\}$ of an $n \times d$-data matrix $\cY$ are mapped onto $\tsd$. 

The empirical covariance matrix reads
\begin{align}
	\hat\Sigma_{\tsd}
	&=\hat\Sigma_{\tsd}(\cY) =\frac{1}{n}\sum_{i=1}^n Log_\mu(y_i)Log_\mu(y_i)^T. 
	\label{eq:empirical_covariance_matrix_ts}
\end{align}
In analogy to the Mahalanobis transformation in $\Rd$ the transformed vector in $\tsd$ reads
\begin{align}
	v_i 
	&= 
	\hat\Sigma_{\tsd}(\cY)^{-1/2}Log_\mu(y_i).         
	\label{eq:mahala-ts}
\end{align}
Note that the condition $||v_i||_2 <\pi$, required for the application of the exponential map, may not be fulfilled. To ensure this, we normalise $\hat\Sigma_{\tsd}(\cY)^{-1/2}$ by
\begin{align}
	\hat\Sigma^*_{\tsd}(\cY)^{-1/2}
	&=\frac{\hat\Sigma_{\tsd}(\cY)^{-1/2}}{||\hat\Sigma_{\tsd}(\cY)^{-1/2}||_2},
\label{eq:hatSigma^*}
\end{align}
where $||\cdot ||_2$ is the spectral norm.
Note that using \eqref{eq:hatSigma^*} in \eqref{eq:mahala-ts} could increase the concentration of the points around $\mu$. 

The Mahalanobis transformation of $y\in\Sd$ reads
\begin{align}
	x 
	&= 
	Exp_\mu \left(\hat\Sigma^*_{\tsd}(\cY)^{-1/2}Log_\mu(y)\right).         
	\label{eq:mahala-trafo-Sd}
\end{align}
For $d=3$, the  Mahalanobis transformation is illustrated in Figure \ref{fig:mahalanobis-trafo-procedure}.

We will use \eqref{eq:mahala-trafo-Sd} to transform realisations $y_1,\dots,y_n\in\Sd$ an elliptically symmetric distribution. If the Mahalanobis-transformed vectors are rotationally symmetric about $\mu$, we will use the results of Ley at al. \cite{Ley14} for the rotationally symmetric case. They are shortly summarised in the following section. 
\subsection{Quantiles for directional data and the angular Mahalanobis depth}\label{sec:Quantiles and depth for directional data}
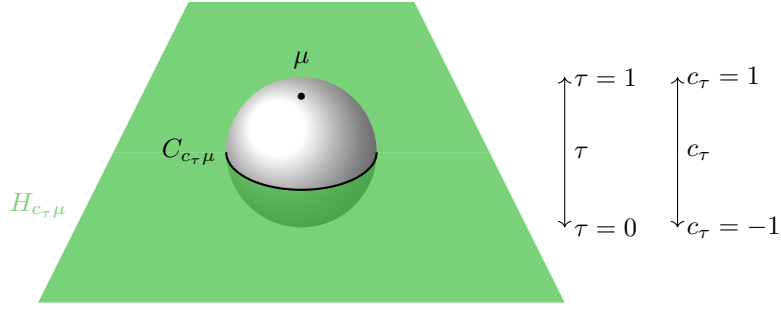
\begin{figure}
	\centering
	\begin{tikzpicture}
		\def\rad{1cm}	
		\fill[green!50!gray, opacity=0.7] (0,0) -- (5,0) -- (4,2) -- (1,2) -- cycle;
		\shade[ball color=white] (2.5,0) circle (1);
		\fill[green!50!gray, opacity=0.7] (0,0) -- (1.5,0) arc (180:360:1 and 0.5) -- (5,0) -- (6,-2) -- (-1,-2) -- cycle;
        \draw[thick]
          (3.5,0) arc [start angle=0,end angle=-180,x radius=\rad,y radius=5mm]node[left] {$C_{c_\tau \mu}$};		
		\draw (2.5,1) node[above] {$\mu$};
		\draw (2.5,0.5) node[above] {\textbullet};
		\draw[green!50!gray, opacity=0.7] (-1,-1) node[above] {$H_{c_\tau\mu}$};
		\draw [->] (6,0) -- (6,1) node[right] {$\tau=1$};
		\draw [<-] (6,-1) -- (6,0) node[right] {$\tau$};
		\draw (6,-1) node[right] {$\tau=0$};
		\draw [->] (7.5,0) -- (7.5,1) node[right] {$c_\tau=1$};
		\draw [<-] (7.5,-1) -- (7.5,0) node[right] {$c_\tau$};
		\draw (7.5,-1) node[right] {$c_\tau=-1$};
	\end{tikzpicture}
	\caption{The hyperplane $H_{c_\tau\mu}$ divides the unit sphere $\S$ into 
	an upper and a lower quantile cap.
		The intersection of $H_{c_\tau\mu}$ with $\S$ corresponds to the $\tau$-depth contour denoted by $C_{c_\tau \mu}$.}
	\label{fig:quantile_projection}
\end{figure}
The concept of quantiles for directional data and the angular Mahalanobis depth from Ley et al. \cite{Ley14} are summarised in the following. 
\subsubsection{Quantiles for directional data}
The quantile check function, known from quantile regression \cite{Koenker05}, reads 
$\rho_\tau(z) = z(\tau - \one[z\leq 0])$, where $z\in\mathbb{R}$, $\tau\in[0,1]$ and $\one[z\leq 0]$ the indicator function.
The projection quantile, i.e., the univariate $\tau$-quantile of $X^T\mu$, is defined by \cite{Ley14,Kong12}
\begin{align}
	c_\tau &= \arg\underset{c\in[-1,1]}{\min} E[\rho_\tau(X^T\mu - c)].  
	\label{eq:quantile_projection}
\end{align}
The partition of the sphere induced by the hyperplane 
\begin{align}
	H_{c_\tau\mu}=\{ x\in\Rd | c_\tau=x^T\mu\}   
	\label{eq:hyper-plane-H} 
\end{align}
defines the $\tau$-depth contour 
\begin{align}
	C_{c_\tau \mu}=H_{c_\tau\mu}\cap\Sd. 
	\label{eq:depth-contour}  
\end{align}
Figure \ref{fig:quantile_projection} illustrates $H_{c_\tau\mu}$ and $C_{c_\tau \mu}$ for $d=3$.
Note that the $0.5$-quantile $c_{0.5}$ is not related to the most central point.
Its associated hyperplane $H_{c_{0.5}\mu}$ divides the probability mass into two equal halves.
The empirical projection quantile \cite{Ley14} reads
\begin{align}
	\hat{c}_\tau = \arg\underset{c\in[-1,1]}{\min} \sum_{i=1}^N[\rho_\tau(X_i^T\hat{\mu} - c)].
	\label{eq:empirical_proj_quant}
\end{align}

\paragraph{Asymptotic properties}
The following Bahadur-type representation of $\hat{c}_\tau$ is proven in \cite{Ley14} (with slightly differing notation).

\begin{proposition}[Proposition 3.1 in \cite{Ley14}]\label{prop:Ley_prop_3_1}
	Let $F\in \cF$ and let $f_{proj}$ denote the common density of the projections $X^T_i\mu$, $i=1,\dots,n$. Set $\Delta_{c_\tau}:=f_{proj}(c_\tau)$. 
	Then there exists a $d$-vector $\Delta_{\mu,c_\tau}$ such that
	\begin{align}
		n^{1/2}(\hat{c}_\tau - c_\tau) 
		= 
		\frac{n^{1/2}}{\Delta_{c_\tau}} \sum_{i=1}^N (\tau - \one[X^T_i\mu\leq c_\tau])] 
		- \frac{\Delta_{\mu,c_\tau}^T}{\Delta_{c_\tau}} n^{1/2} (\hat{\mu}-\mu)+o_P(1)
		\label{eq:Bahadur}
	\end{align}
	as $n\rightarrow\infty$ under the joint distribution of $X_1,\dots,X_n.$
\end{proposition}
Note that by \eqref{eq:Bahadur}, the rather complicated non-linear estimator $\hat{c}_\tau$ can be represented as a sum of i.i.d. random variables and the scaled difference between $\mu$ and its estimator $\hat{\mu}$. 
However, the calculation of the $d$-vector $\Delta_{\mu,c_\tau}=\frac{d}{dc}E\left((\tau-1[X_i^T\mu\leq c])X_i\right)_{|c=c_\tau}$ (see the proof of Proposition 3.1 in \cite{Ley14}) is not straightforward.
In the rotationally symmetric case, the representation in Equation (\ref{eq:Bahadur}) simplifies.

\begin{proposition}[Proposition 3.2 in \cite{Ley14}]\label{prop:Ley-rot-sym}
	Let $F\in \cR$.
	Then 
	\begin{equation}
		n^{1/2}(\hat{c}_\tau - c_\tau) 
		= 
		\frac{n^{1/2}}{\Delta_{c_\tau}} \sum_{i=1}^N (\tau - \one[X^T_i\mu\leq c_\tau])] 
		+o_P(1)
		\label{eq:Bahadur_rot}
	\end{equation}
	as $n\rightarrow\infty$ under the joint distribution of $X_1,\dots,X_n.$
	Therefore, letting $f_{proj}$ stand for the density of $X_i^T\mu$, we have that $n^{1/2}(\hat{c}_\tau-c_\tau)$ is asymptotically normal with mean zero and variance $\frac{(1-\tau)\tau}{f^2_{proj}(c_\tau)}$.
\end{proposition}

The reason for the simplification in Equation (\ref{eq:Bahadur}) is that $\Delta_{\mu,c_\tau}^T n^{1/2} (\hat{\mu}-\mu) \in o_P(1)$ for $F\in\cR$.
The absence of $\hat\mu$ in Equation (\ref{eq:Bahadur_rot}) means that any root-$n$ consistent estimator (e.g., the empirical Fisher spherical median $\hat{\mu}$ or the spherical mean $\sum_{i=1}^{n}X_i/||\sum_{i=1}^{n}X_i||_2$) can substitute $\mu$ without changing the asymptotic distribution, independently of the dimension $d$.
Furthermore, (\ref{eq:Bahadur}) is a Bahadur-type representation for univariate sample quantiles \cite{Bahadur66}. Hence, the directional quantiles of \cite{Ley14} have similar asymptotic properties as the quantiles in $\mathbb{R}$.
Therefore, the directional quantiles of \cite{Ley14} can be regarded as canonical for $F\in\cR$.

\subsubsection{The angular Mahalanobis depth}\label{sec:The angular Mahalanobis depth}
The angular Mahalanobis depth (AMHD) is defined by \cite{Ley14}
\begin{equation}
	AMHD_F(x) = \frac{1}{1+\sfrac{1}{D_F(x)}}  \in \left[0,\sfrac{1}{2}\right]
\end{equation}
where 
\begin{align}
	D_F(x) = \arg\underset{\tau\in[0,1]}{\min}\{c_\tau\geq x^T\mu\}.
	\label{eq:D_F}
\end{align}

It provides a centre-outward ordering by assigning each $x\in \Sd$ its depth value.
The angular Mahalanobis depth is leaned on the classical Mahalanobis depth 
\begin{equation}
	MHD_F (x) = \frac{1}{1+(x-\mu(F))^T(\Sigma(F))^{-1}(x-\mu (F))} , \quad  x\in\Rd.
	\label{eq:Classical_Mahalanobis_depth}
\end{equation}

$\mu(F)$ and $\Sigma(F)$ are location and scatter functionals under $F$, respectively.
The spherical centre $\mu$ corresponds to the centre $\mu(F)$.
$MHD_F$ is suited for elliptically symmetric distributions on $\Rd$ since $\Sigma(F)$ contains all necessary information about the principal axes.
In contrast, $AMHD_F$ is not suited for distributions which are elliptically symmetric about $\mu\in\Sd$ since information about the shape of the distribution is lost due to the projection $X^T\mu$ in the definition of $D_F(x)$.

\section{Quantiles for directional data from elliptically symmetric distributions and the elliptical Mahalanobis depth}\label{sec:Quantiles-ellipse}
Here, we present canonical quantiles and a depth for directional distributions which are elliptically symmetric about $\mu\in\Sd$.
The idea is to transform the elliptical contours in the tangent space to rotationally symmetric contours analogously to \eqref{eq:Elliptical symmetry is connected to rotational symmetry}, such that we are again in the canonical case of Ley et al. \cite{Ley14}.

In the following, we consider random vectors $Y,Y_1,\dots,Y_n\in \Sd$ i.i.d. following a distribution which is elliptically symmetric about $\mu$. Let
\begin{align}
	(\Sigma^*)^{-1/2}=\frac{\Sigma^{-1/2}}{||\Sigma^{-1/2}||_2}
	\label{eq:(Sigma^*)^-1/2}
\end{align}
and
\begin{align}
	(\Sigma^*)^{1/2}=||\Sigma^{-1/2}||_2\cdot \Sigma^{1/2}.
	\label{eq:(Sigma^*)^1/2}
\end{align}
Then, $||(\Sigma^*)^{-1/2}Log_\mu(Y)||_2  < \pi$ such that we can define
\begin{align}
	\cG(Y) = Exp_\mu\left((\Sigma^*)^{-1/2}Log_\mu(Y)\right).
	\label{eq:cG(Y)}
\end{align}
$\cG(Y)$ is locally diffeomorphic since $Exp_\mu$ and $Log_\mu$ are locally diffeomorphic, and $(\Sigma^*)^{-1/2}$ is invertible. Its inverse is
\begin{align}
	\cG^{-1}(Y) = Exp_\mu\left((\Sigma^*)^{1/2}Log_\mu(Y)\right).
	\label{eq:cG^-1(Y)}
\end{align}
\subsection{Quantiles for directional data from elliptically symmetric distributions}
\begin{figure}
    \centering
    \includegraphics{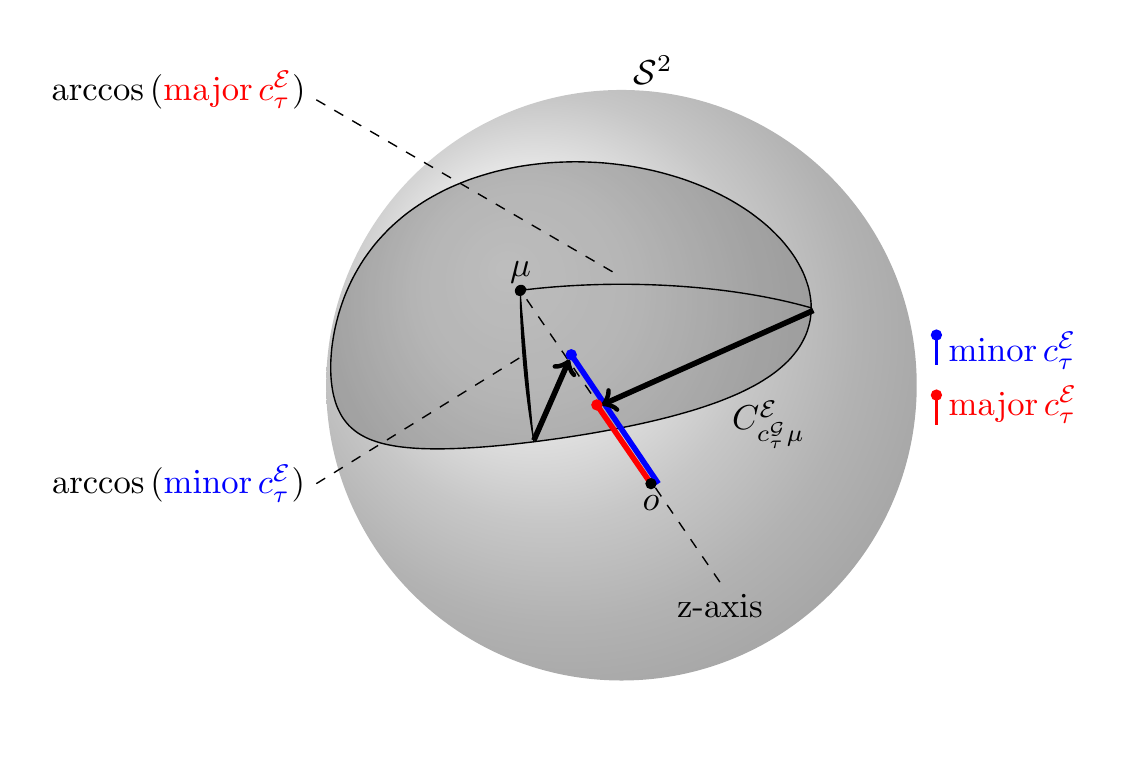}
    \caption{Illustration of $\minor c^\cE_\tau$, $\major c^\cE_\tau$ and the intrinsic small semi-axis $\arccos\left(\minor c^\cE_\tau\right)$ and the intrinsic large semi-axis $\arccos\left(\major c^\cE_\tau\right)$ of $C^\cE_{c^\cG_\tau\mu}$ for some $\tau$ and $d=3$.}
    \label{fig:Sphere_with_elliptic_contour}
\end{figure}
We define the elliptical projection quantile by
\begin{align}
	c^\cG_\tau &= \arg\underset{c\in[-1,1]}{\min} E[\rho_\tau(\cG(Y)^T\mu - c)].
	\label{eq:c^cG_tau}
\end{align}
The partition of the sphere induced by the hyperplane $H_{c^\cG_\tau\mu}=\{ x\in\Rd | c^\cG_\tau=\cG(Y)^T\mu\}$ defines the $\tau$-depth contour $C_{c^\cG_\tau\mu}$ as in \eqref{eq:depth-contour}.
Note that $c^\cG_\tau = c_\tau$ if $\Sigma^*=I_{d-1}$ since then $\cG(Y) =Y$.

An elliptical depth contour is obtained from $C_{c^\cG_\tau\mu}$ by inverting the transformation via the tangent space shown in Figure \ref{fig:mahalanobis-trafo-procedure}, that is 
\begin{align}
	C^\cE_{c^\cG_\tau\mu}
	&= 
	\left\{\cG^{-1}(Y) | Y\in C_{c^\cG_\tau\mu} \right\}.
	\label{eq:elliptical-countour-depth}
\end{align}
Note that $C^\cE_{c^\cG_\tau\mu}$ can be elliptically shaped which is not the case for $C_{c^\cG_\tau\mu}$.

To define an equivalent to the semi-axes lengths of an ellipse, we set
\begin{equation}
	\minor c^\cE_\tau
	=\max_{x\in C^\cE_{c^\cG_\tau\mu}} x^T\mu	\text{ and }
	\major c^\cE_\tau
	=\min_{x\in C^{\cE}_{c^{\cG}_\tau\mu}} x^T\mu. 	
	\label{eq:minormajor}
\end{equation}
The intrinsic semi-minor axis of $C^\cE_{c^\cG_\tau\mu}$ is $\arccos\left(\minor c^\cE_\tau\right)$ and the intrinsic semi-major axis of $C^\cE_{c^\cG_\tau\mu}$ is $\arccos\left(\major c^\cE_\tau\right)$.
$\minor c^\cE_\tau$ and $\major c^\cE_\tau$ contain the main information about the concentration and shape of the distribution of $Y$ around $\mu$. 
A large difference between $\minor c^\cE_\tau$ and $\major c^\cE_\tau$ indicates a strong deviation from a rotationally symmetric distribution, 
whereas $\minor c^\cE_\tau = \major c^\cE_\tau$ in the rotationally symmetric case.
See Figure \ref{fig:Sphere_with_elliptic_contour} for an illustration for $d=3$.
Note that $\major c^\cE_\tau \leq c_\tau \leq \minor c^\cE_\tau$ by construction.

The empirical elliptical projection quantile reads
\begin{align}
	\hat{c}^\cG_\tau = \arg\underset{c\in[-1,1]}{\min} \sum_{i=1}^N[\rho_\tau(\hat\cG(Y_i)^T\hat{\mu} - c)]
	\label{eq:empirical_proj_quant_ell}
\end{align}
with 
\begin{align}
	\hat\cG(Y) 
	&= Exp_\mu((\hat\Sigma^*)^{-1/2}(Log_\mu(Y)))
\end{align}
and $(\hat\Sigma^*)^{-1/2}$ given in \eqref{eq:hatSigma^*}. The empirical versions of \eqref{eq:elliptical-countour-depth} and \eqref{eq:minormajor} are denoted by $C^\cE_{\hat c^\cG_\tau\mu}$, $\minor \hat{c}^\cE_\tau$ and $\major \hat{c}^\cE_\tau$, respectively.

\subsection{The elliptical Mahalanobis depth}
The elliptical Mahalanobis depth (EMHD) is defined by
\begin{equation}
	EMHD_{F}(y) = \frac{1}{1+\sfrac{1}{D^\cG_{F}(y)}} 
	= \frac{D^\cG_{F}(y)}{1+D^\cG_{F}(y)} 
	\in \left[0,\sfrac{1}{2}\right],    
\end{equation}
where 
\begin{align}
	D^\cG_{F}(y) = \arg\underset{\tau\in[0,1]}{\min}\{c^\cG_\tau\geq \cG(y)^T\mu\}.
\end{align}
As the angular Mahalanobis depth, the elliptical Mahalanobis depth is leaned on the Classical Mahalanobis depth $MHD_F$ given in Equation (\ref{eq:Classical_Mahalanobis_depth}).
Our approach additionally yields a connection to $\Sigma(F)$ via $\hat\Sigma^*$ which contains all necessary information about the principal axes.
Furthermore, $EMHD_F$ contains $AMHD_F$ as a special case: They are equal if the depth contours are rotationally symmetric.

\section{Applications}\label{seq:Applications Quantiles}
To confirm our findings, we perform a Monte Carlo simulation study. 
We generated 
four independent samples 
\begin{align*}
	y_{l,i}, \quad l=1,2,3,4, \quad i=1,\dots,200,
\end{align*}
of Kent distributions with $\mu=(0,0,1)^T$, $A=\diag(\beta,-\beta,0)$, and 
$\kappa=5$, $\beta=2$ ($l=1$),
$\kappa=7$, $\beta=3$ ($l=2$), 
$\kappa=10$, $\beta=4$ ($l=3$), and
$\kappa=12$, $\beta=5$ ($l=4$).

%
The Mahalanobis-transformed $y_{l,i}$ are denoted by $x_{l,i}$. The longitude of $y_{l,i}$ is denoted by  $\phi_{y_{l,i}}$, and the longitudes of $x_{l,i}$ are $\phi_{x_{l,i}}$.

\begin{figure}
	\centering
	\begin{subfigure}{0.24\linewidth}
		\includegraphics[width=\linewidth]{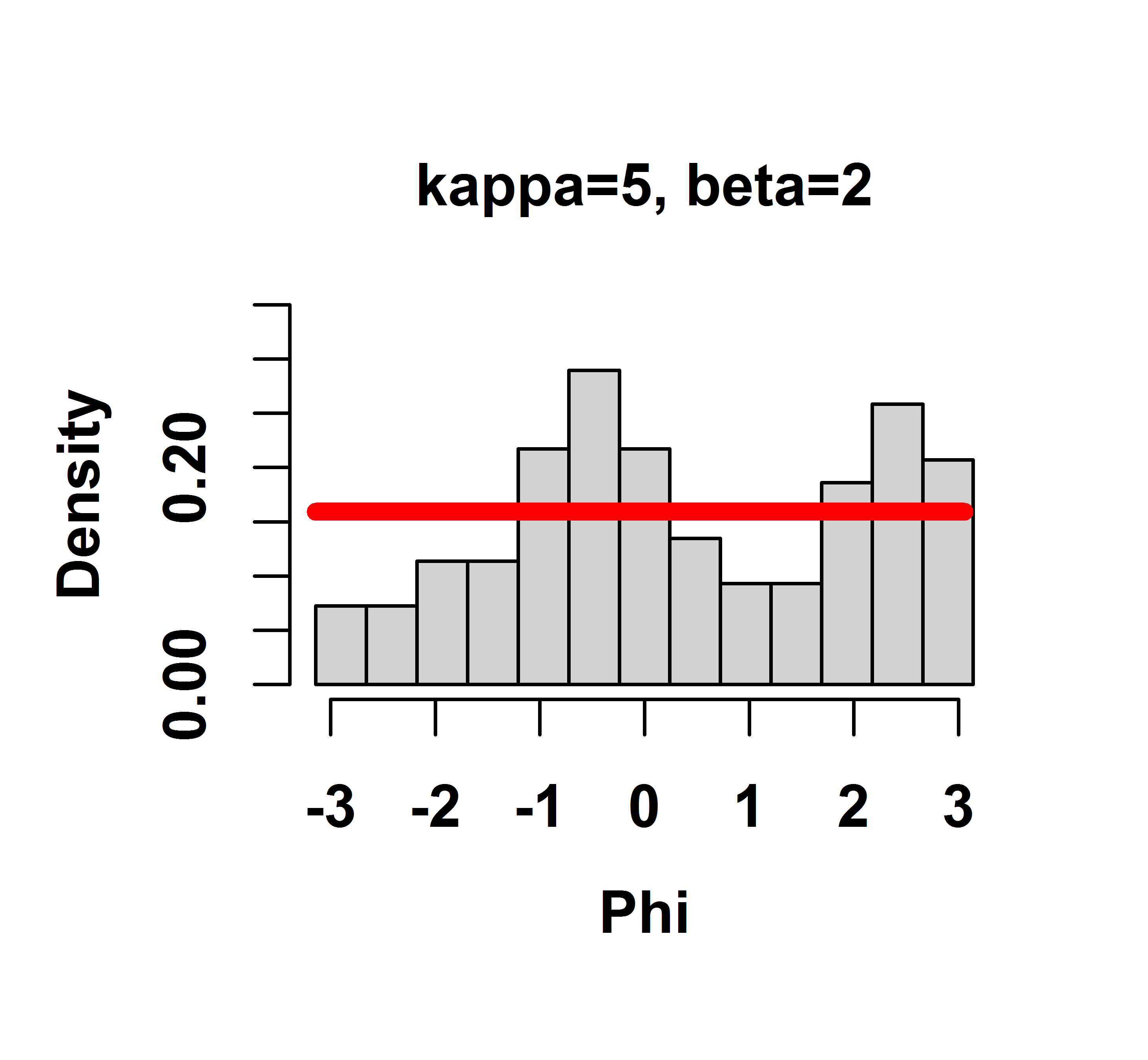}
		\caption{$\phi_{y_{1,i}}$}
	\end{subfigure}
	\begin{subfigure}{0.24\linewidth}
		\includegraphics[width=\linewidth]{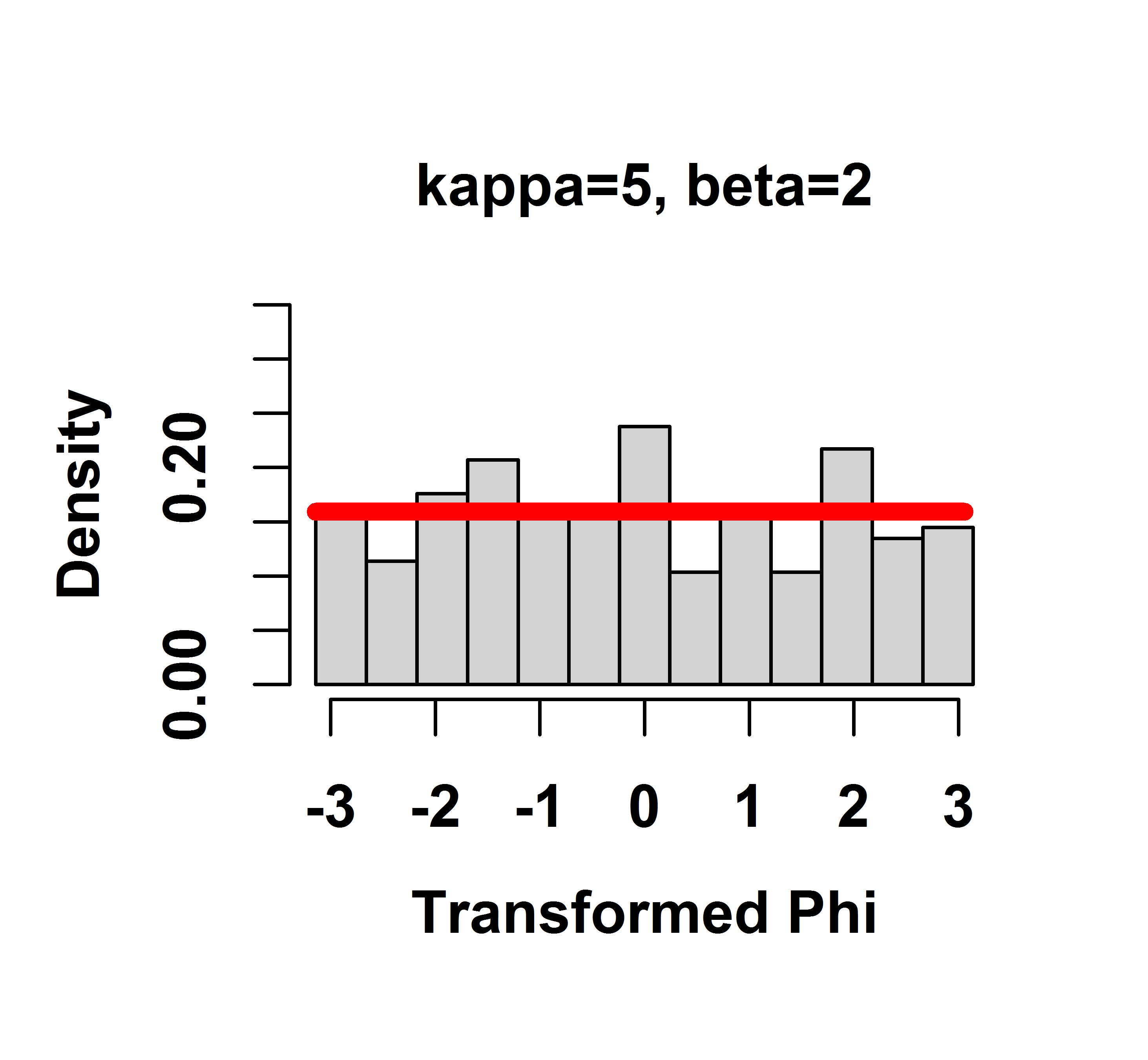}
		\caption{$\phi_{x_{1,i}}$}
	\end{subfigure}
	\begin{subfigure}{0.24\linewidth}
		\includegraphics[width=\linewidth]{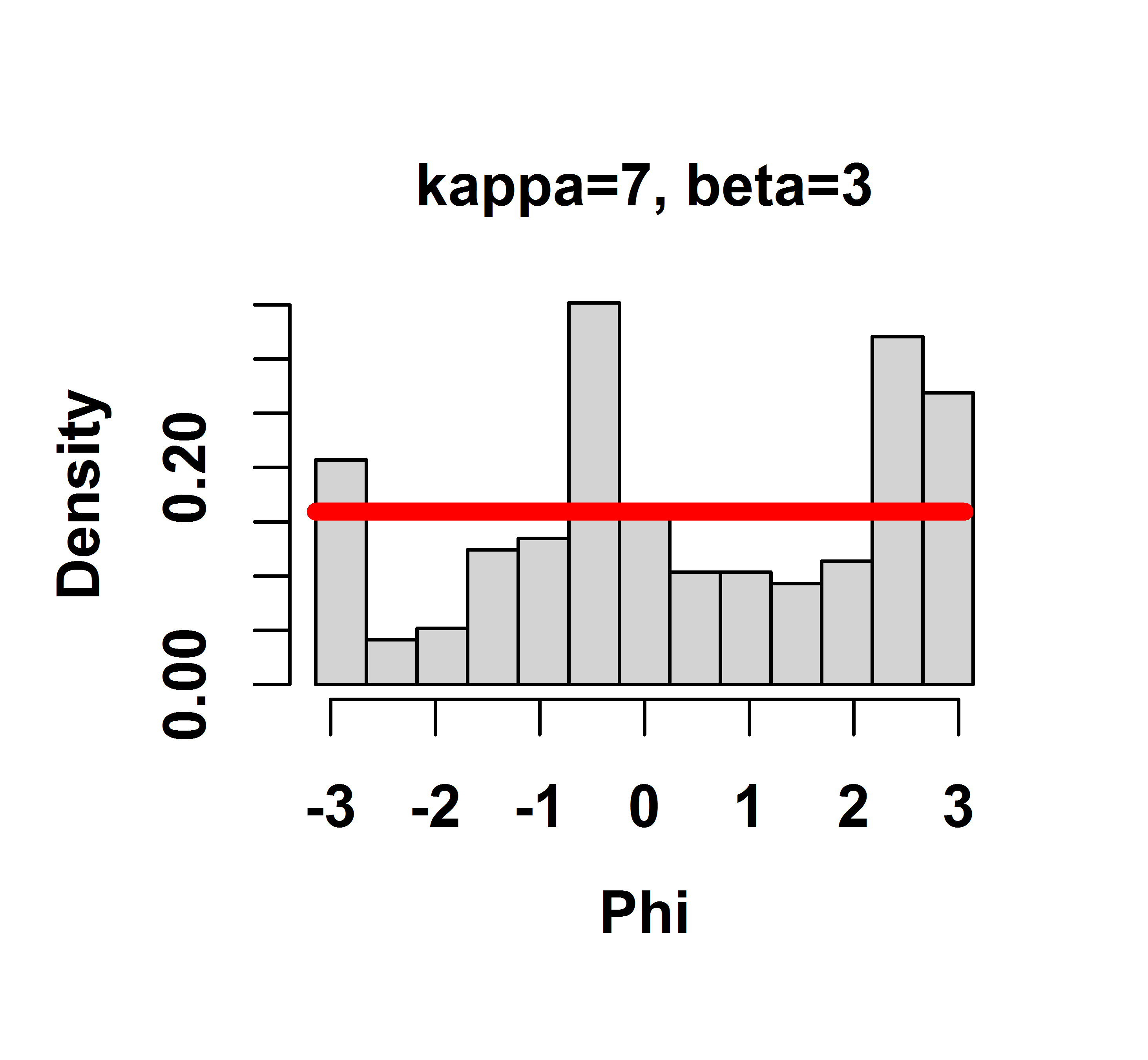}
		\caption{$\phi_{y_{2,i}}$}
	\end{subfigure}
	\begin{subfigure}{0.24\linewidth}
		\includegraphics[width=\linewidth]{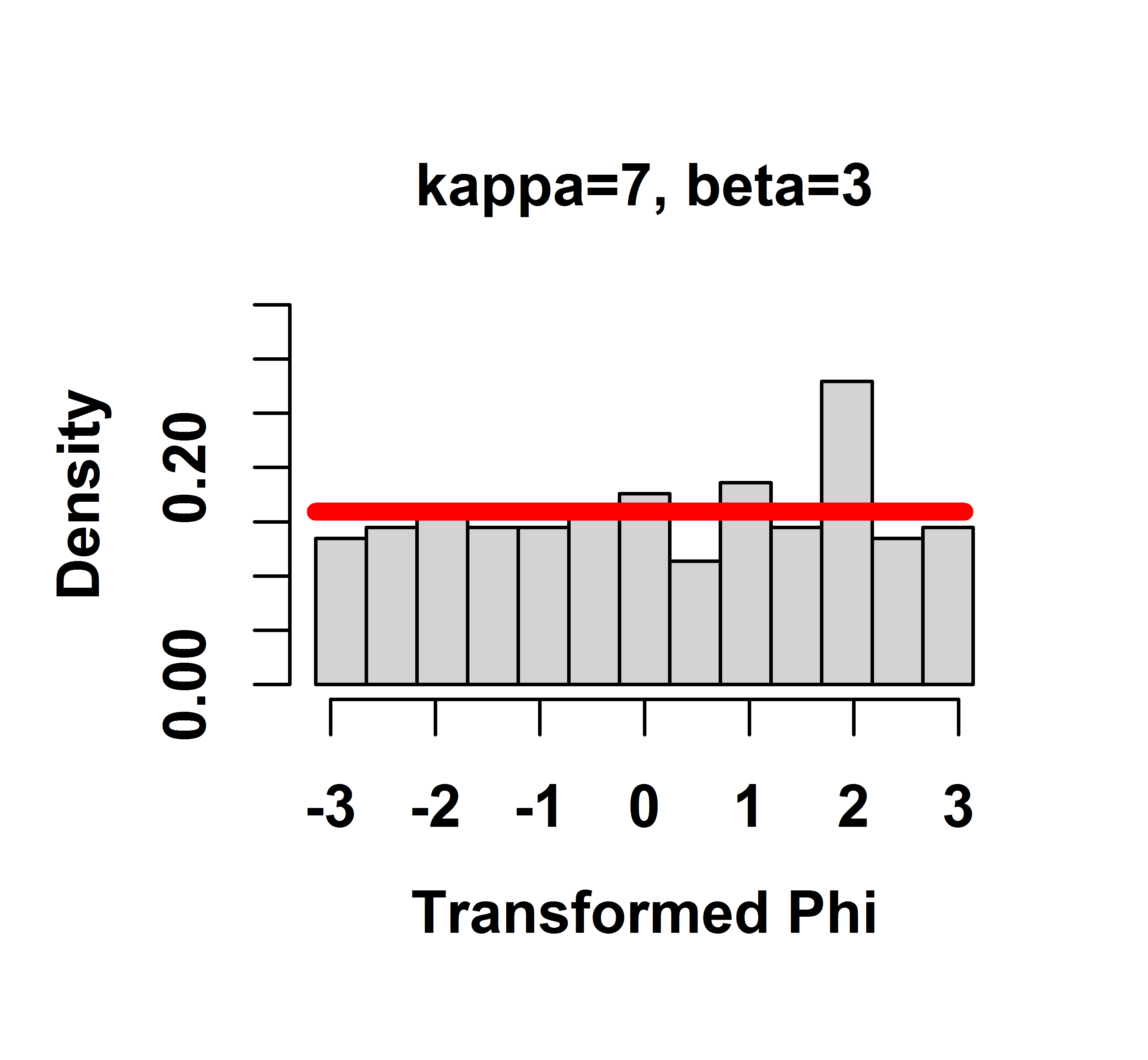}
		\caption{$\phi_{x_{2,i}}$}
	\end{subfigure}
	\begin{subfigure}{0.24\linewidth}
		\includegraphics[width=\linewidth]{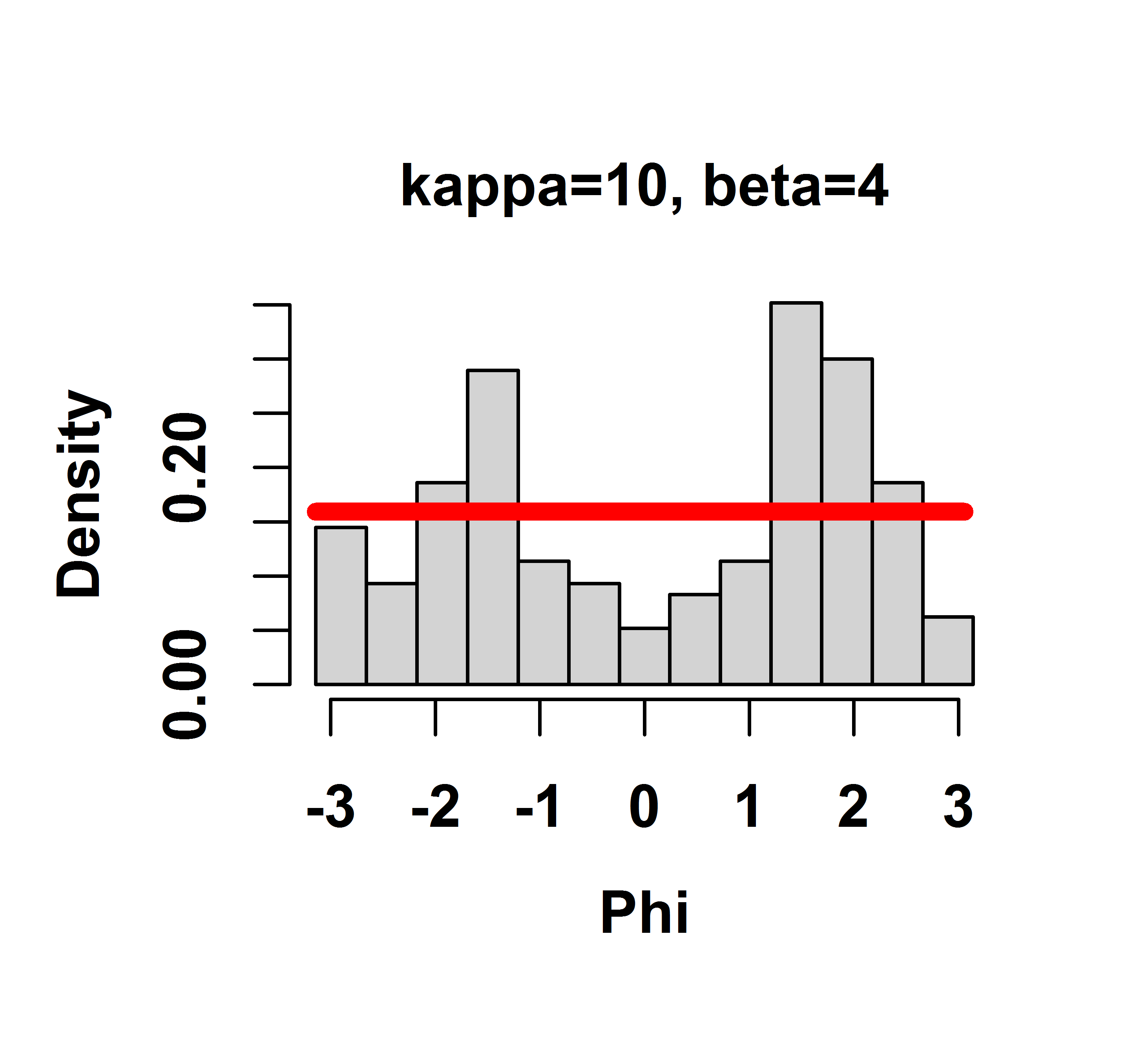}
		\caption{$\phi_{y_{3,i}}$}
	\end{subfigure}
	\begin{subfigure}{0.24\linewidth}
		\includegraphics[width=\linewidth]{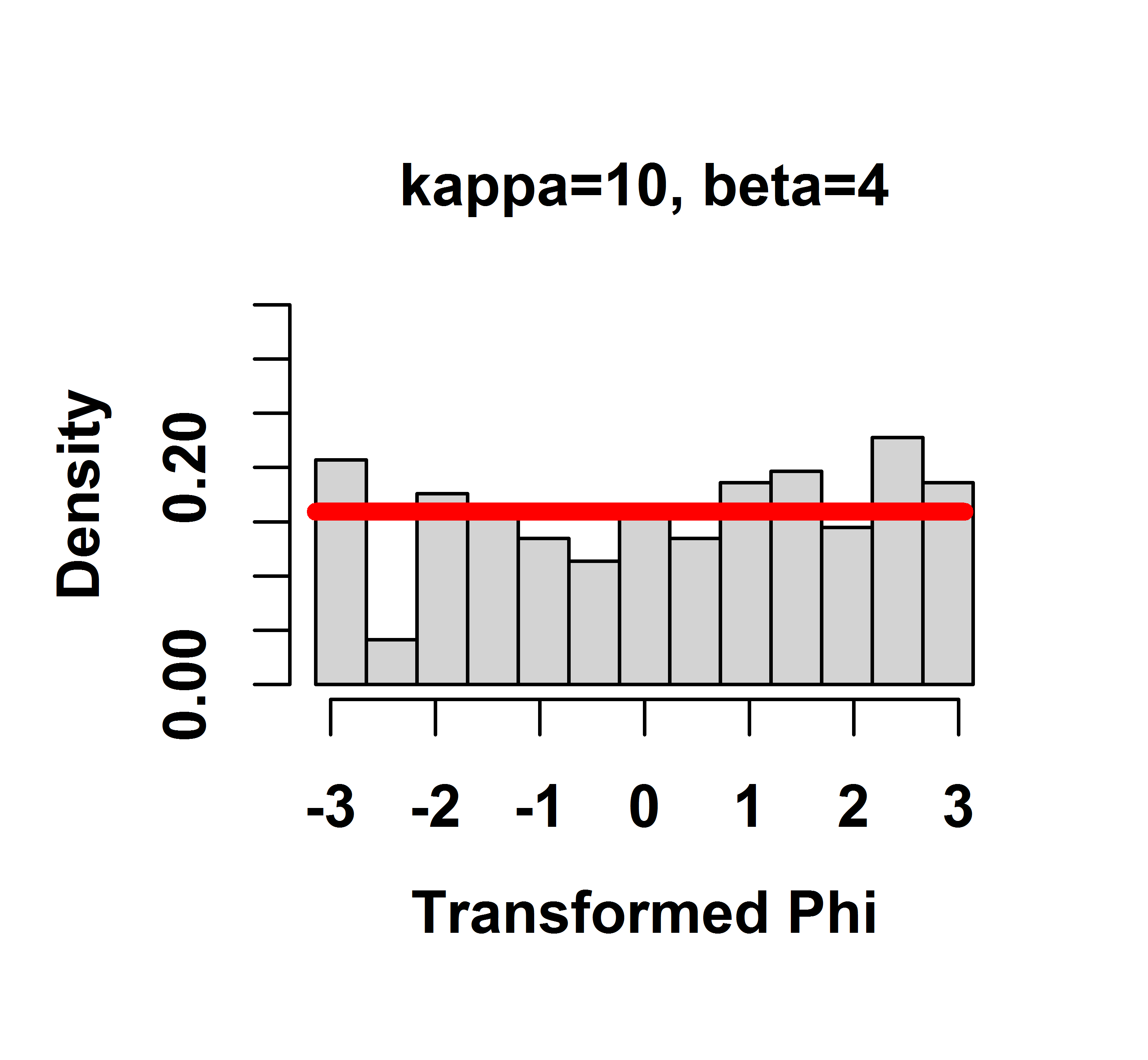}
		\caption{$\phi_{x_{3,i}}$}
	\end{subfigure}
	\begin{subfigure}{0.24\linewidth}
		\includegraphics[width=\linewidth]{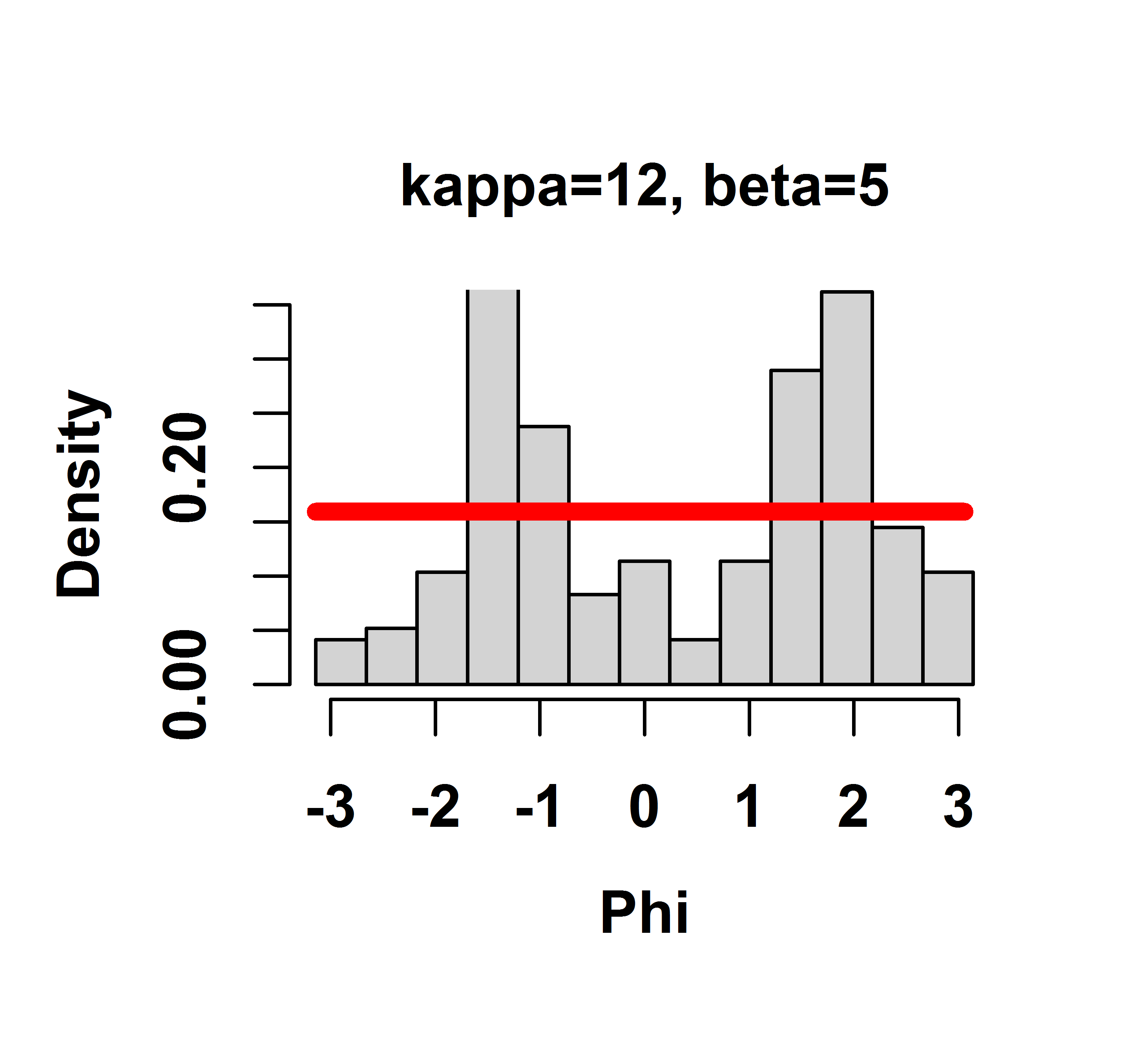}
		\caption{$\phi_{y_{4,i}}$}
	\end{subfigure}
	\begin{subfigure}{0.24\linewidth}
		\includegraphics[width=\linewidth]{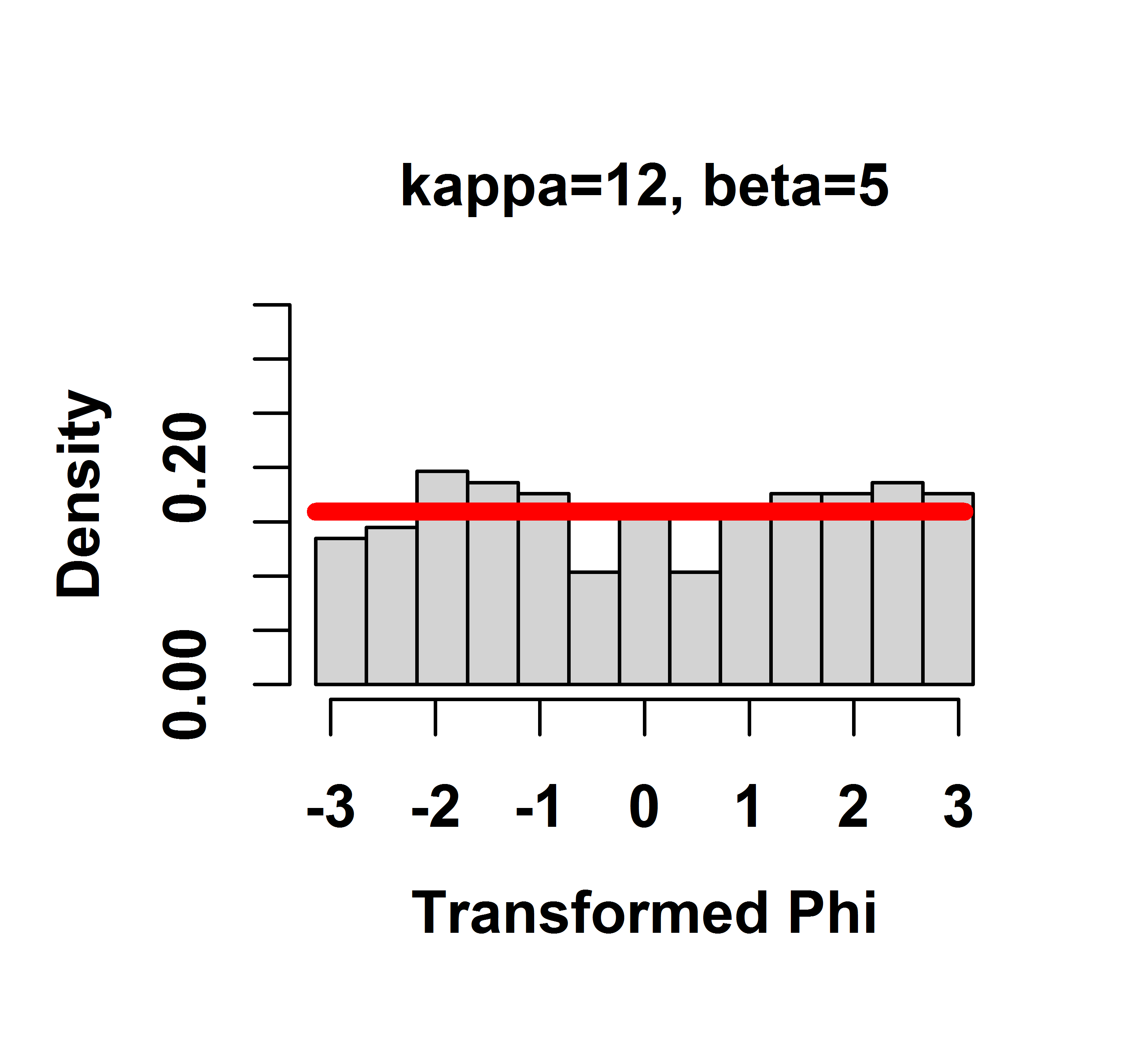}
		\caption{$\phi_{x_{4,i}}$}
	\end{subfigure}
	\caption{
		Histograms of the longitudes $\phi_{y_{l,i}}$ and $\phi_{x_{l,i}}$, $l=1,2,3,4$.
		The red line corresponds to the density of the uniform distribution on $[-\pi, \pi]$.
	}
	\label{fig:histogram_phi_before_after_mahala}
\end{figure}

The histograms of the longitudes shown in Figure \ref{fig:histogram_phi_before_after_mahala} indicate that the transformation leads to uniformly distributed longitudes. This supports that the $x_{l,i}$ are rotationally symmetric about $\mu$.
To confirm this visual impression, we test the hypothesis of uniform longitudes.
We use Watson's test \cite[p. 156]{Jammalamadaka01} implemented in the $R$-package $Directional$ \cite{Directional}.
Watson's test applied on $\phi_{y_{i}}$ gave p-values less than 0.004 for all designs $l=1,2,3,4$.
The p-values for the Mahalanobis-transformed angles $\phi_{x_{i}}$ were 0.7086 ($l=1$), 0.5132 ($l=2$), 0.3436 ($l=3$), 0.5268 ($l=4$) which supports an assumption of uniform longitudes.

\begin{figure}
	\centering
	\begin{subfigure}{0.32\linewidth}
		\includegraphics[width=\linewidth]{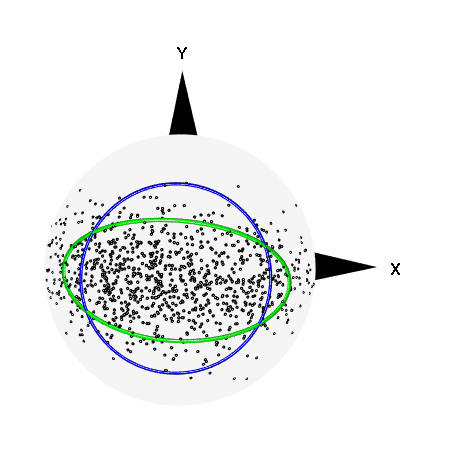}
		\caption{$\tau=0.25$}
	\end{subfigure}
	\begin{subfigure}{0.32\linewidth}
		\includegraphics[width=\linewidth]{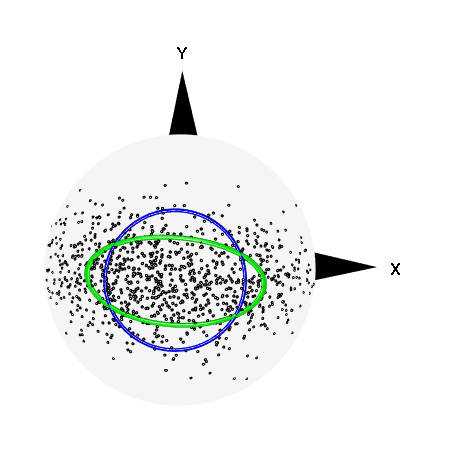}\caption{$\tau=0.5$}
	\end{subfigure}
	\begin{subfigure}{0.32\linewidth}
		\includegraphics[width=\linewidth]{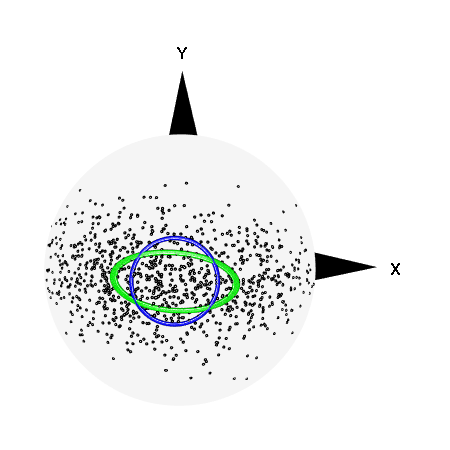}\caption{$\tau=0.75$}
	\end{subfigure}
	\caption{
		Realisations $y_{4,i}$, $i=1,\dots,n$, are given as points on $\S$ with empirical $\tau$-depth contours, $\tau=0.25,0.5,0.75$.
		The blue circle corresponds to the empirical $\tau$-depth contour $C_{\hat c_\tau \mu}$.
		The green ellipse corresponds to the empirical $\tau$-depth contour $C^\cE_{\hat c^\cG_\tau\mu}$.
		The Z-direction points out of the page.
	}
	\label{fig:quantiles_comparsion}
\end{figure}

For illustration, the empirical quartiles $\hat c_\tau$, $\tau=0.25, 0.5,$ and $0.75$ as well as $\minor \hat c^\cE_\tau$ and $\major \hat c^\cE_\tau$ for $l=4$ are given in Table \ref{tab:quantiles_comparsion}.  
Figure \ref{fig:quantiles_comparsion} shows the corresponding depth contours. 

\begin{table}
	\centering
	\begin{tabular}{c|ccc}
		$\tau$ & 0.25 & 0.5 & 0.75 \\ 
		\hline
		$\hat{c}_\tau$          & 0.8110  & 0.9090 & 0.9660  \\
		$\minor\hat{c}^\cE_\tau$  & 0.9370  & 0.9677 & 0.9870  \\
		$\major\hat{c}^\cE_\tau$  & 0.6847  & 0.8347 & 0.9245
	\end{tabular}
	\caption{
		The empirical quartiles $\minor \hat c^\cE_\tau, \major  \hat c^\cE_\tau$ and $\hat c_\tau$ of $y_{4,i}$, $i=1,\dots,n$. The sample is shown in Figure \ref{fig:quantiles_comparsion}.
	}
	\label{tab:quantiles_comparsion}
\end{table}

\subsection{Trimming of directional data}
\begin{figure}
	\centering
	\begin{subfigure}{0.4\linewidth}
		\includegraphics[width=\linewidth]{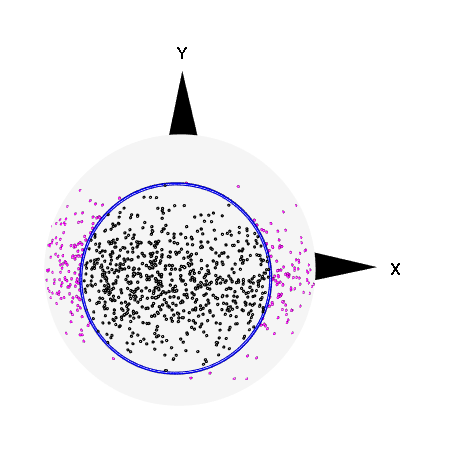}
		\caption{$C_{c_\tau\mu}$}
	\end{subfigure}
	\begin{subfigure}{0.4\linewidth}
		\includegraphics[width=\linewidth]{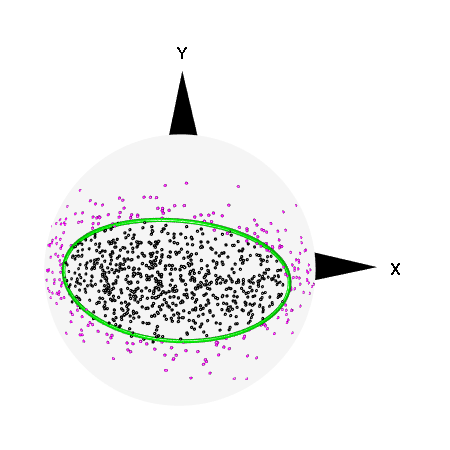}
		\caption{$C^\cE_{c^\cG_\tau\mu}$}
	\end{subfigure}
	\caption{
		Realisations $y_{4,i}$, $i=1,\dots,n$, are given as points on $\S$.
		The blue circle corresponds to $C_{\hat c_\tau\mu}$ (left) and the green ellipse corresponds to $C^\cE_{\hat c^\cG_\tau\mu}$ with $\tau=0.25$ of $y_{4,i}$, $i=1,\dots,n$. 
		Trimmed points are purple. 
		The Z-direction points out of the page.
	}
	\label{fig:trimming}
\end{figure}

The angular Mahalanobis depth $AMHD_F$ is canonical for trimming of directional data from $F\in\cR$.
The trimming corresponds to deleting the points on $\Sd$ below the $\tau$-depth contour $C_{c_\tau \mu}$ given in \eqref{eq:depth-contour} with $\tau\in [0,1]$. 
If the underlying distribution has elliptical contours, trimming results in circular contours when using $AMHD_{F}$. In contrast, trimming based on $EMHD_{F}$, which deletes points below $C^\cE_{c^\cG_\tau\mu}$, preserves the elliptical shape of the contours.
See Figure \ref{fig:trimming} for an illustration.

\subsection{Analysis of fibre directions in ultra-high performance fibre-reinforced concrete}

%
Ultra-high performance fibre-reinforced concrete (UHPFRC) is a relatively new material in civil engineering. 
If cracks appear in the concrete due to loading, fibres crossing the crack counteract the crack propagation. 
As fibres have no directional sense, our data are restricted to the upper hemisphere of $\S$. 
We analyse a data set from \cite{Maryamh20} which consists of $n=598$ measurements of fibre directions.
The fibres crossed a crack in a UHPFRC-specimen subject to a bending test. The crack has a planar shape with normal direction corresponding to the Z-axis used in the analysis.
The fibre directions are denoted by $y_i$, their Mahalanobis transforms are $x_i$, $i=1,\dots,n$.
Furthermore, we denote by $\phi_{y_{i}}$ the longitudes of $y_{i}$ and by $\phi_{x_{i}}$ the longitudes of $x_{i}$.

\subsubsection{Visual inspection, rotational symmetry, and quartiles}
\begin{figure}
	\centering
	\begin{subfigure}{0.48\linewidth}
		\includegraphics[width=\linewidth]{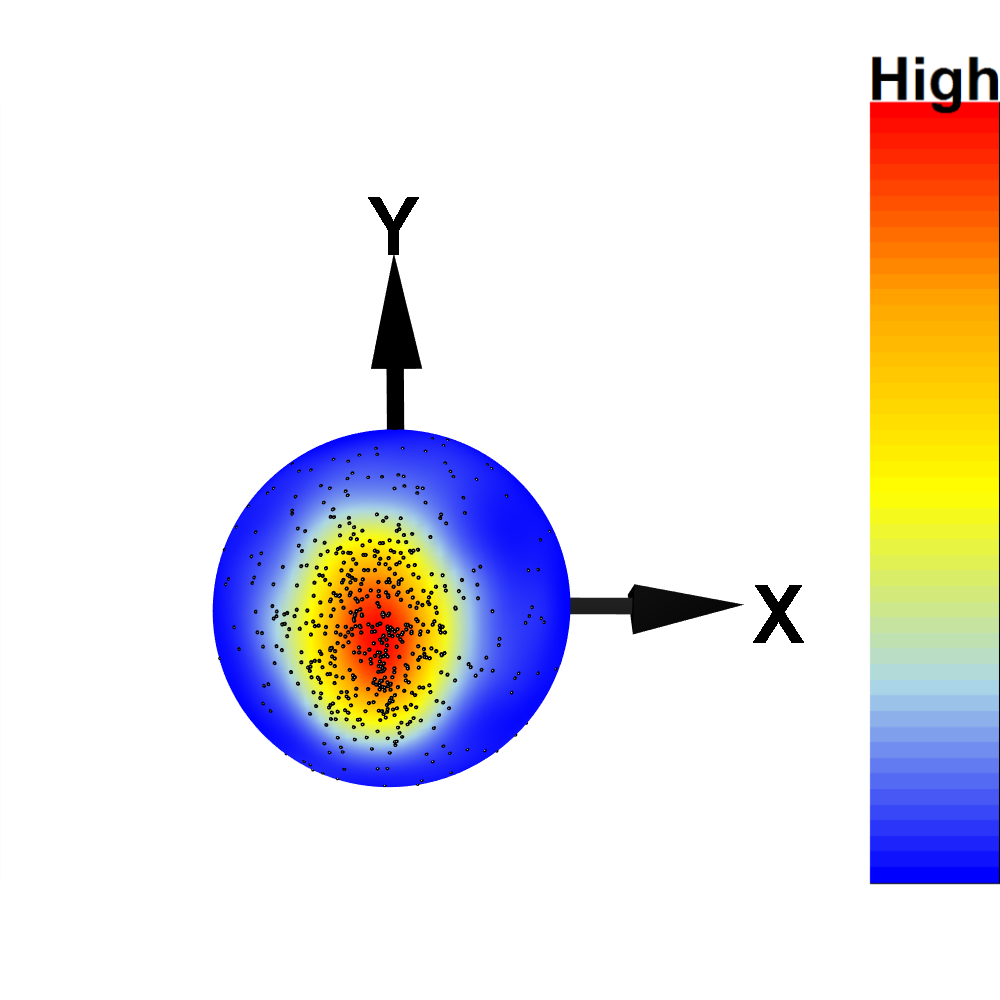}\caption{Original fibre directions $y_i$}\label{fig:UHPFRC-original}
	\end{subfigure}
	\begin{subfigure}{0.5\linewidth}
		\includegraphics[width=\linewidth]{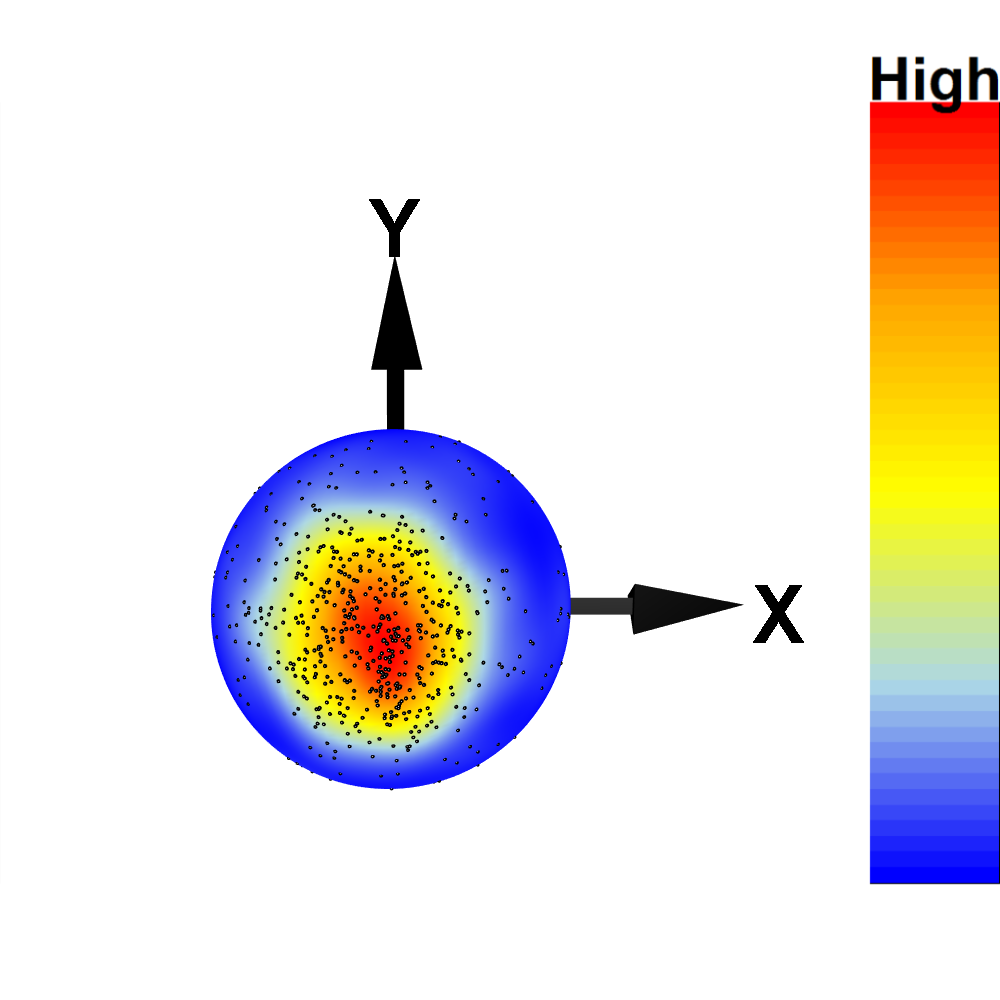}\caption{Mahalanobis transformed directions $x_i$}\label{fig:UHPFRC-mahala}
	\end{subfigure}
	\begin{subfigure}{0.45\linewidth}
		\includegraphics[width=\linewidth]{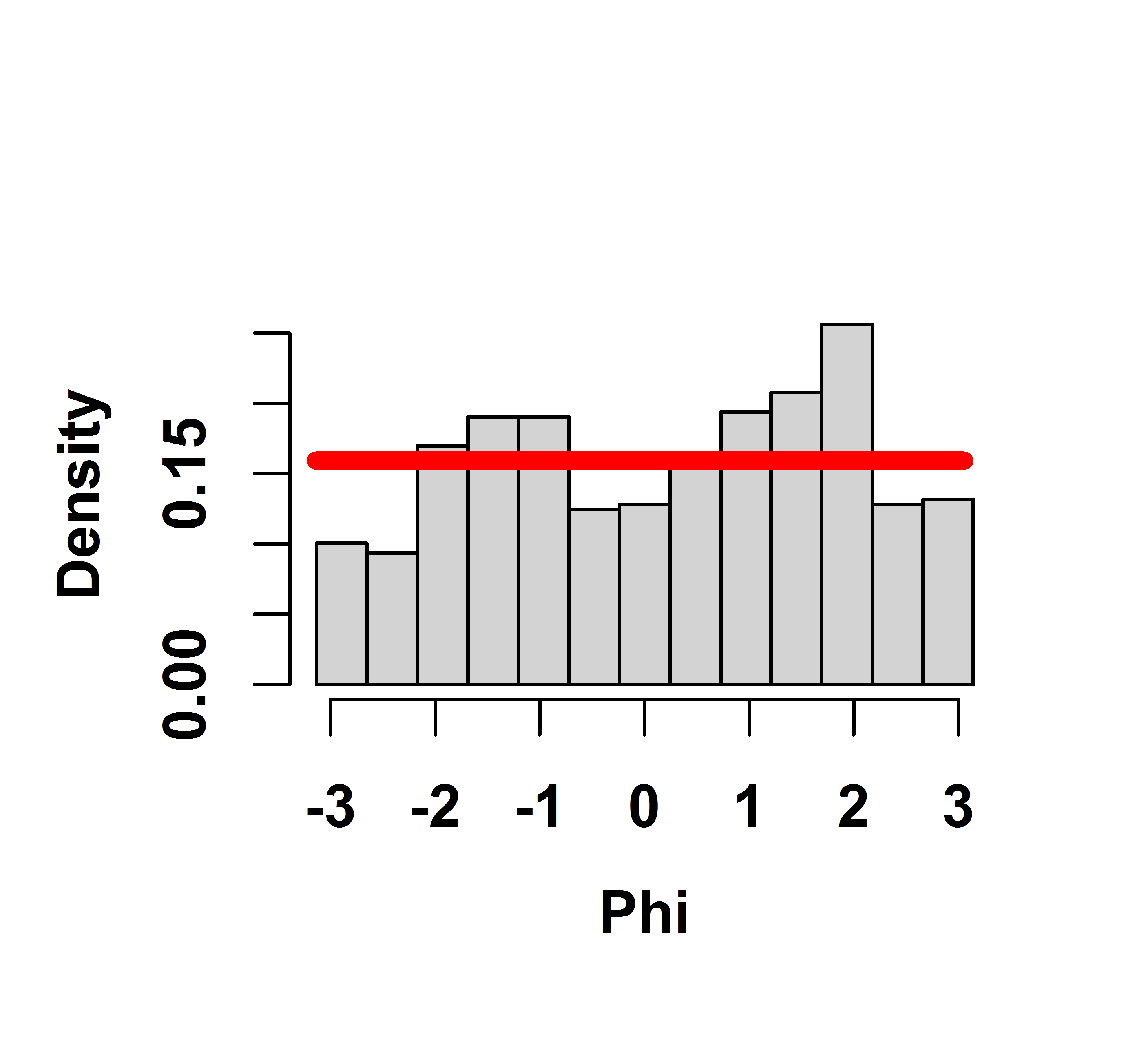}
		\caption{$\phi_{y_{i}}$}\label{fig:UHPFRC-original-phi}
	\end{subfigure}
	\begin{subfigure}{0.45\linewidth}
		\includegraphics[width=\linewidth]{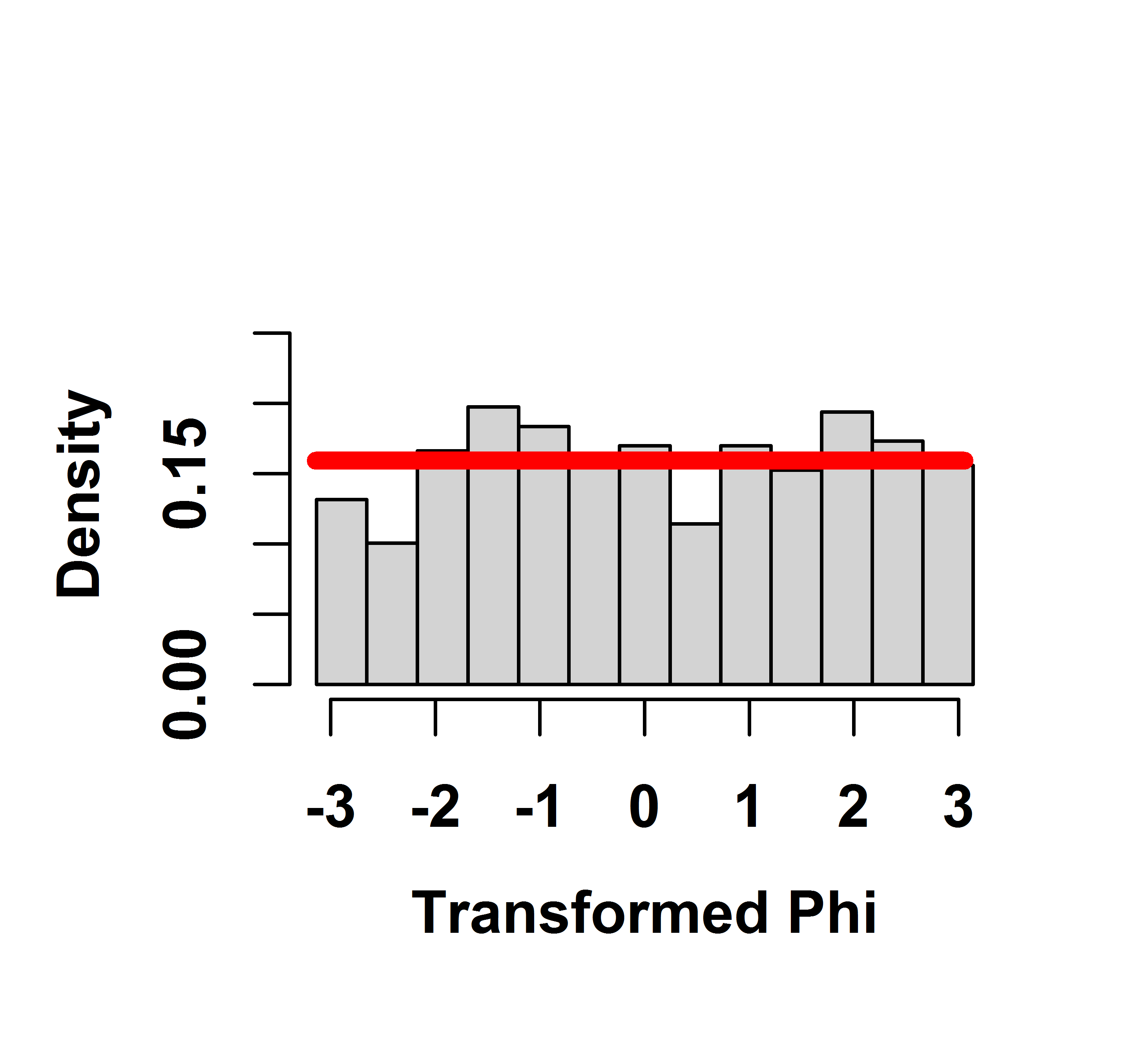}
		\caption{$\phi_{x_{i}}$}\label{fig:UHPFRC-mahala-phi}
	\end{subfigure}
	\caption{
		Fibre directions $y_i$ in UHPFRC before (a) and after (b) Mahalanobis transformation.
		The Z-direction points out of the page.
		Histograms of the longitudes $\phi_{y_{i}}$ (c) and $\phi_{x_{i}}$ (d).
		The red line corresponds to the density of the uniform distribution on $[-\pi, \pi]$.
	}
	\label{fig:UHPFRC}
\end{figure}
\begin{figure}
	\centering
	\begin{subfigure}{0.32\linewidth}
		\includegraphics[width=\linewidth]{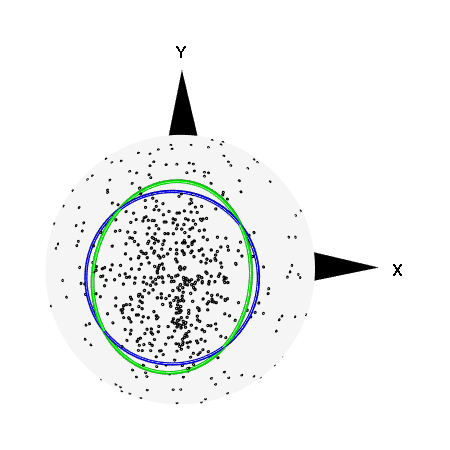}
		\caption{$\tau=0.25$}
	\end{subfigure}
	\begin{subfigure}{0.32\linewidth}
		\includegraphics[width=\linewidth]{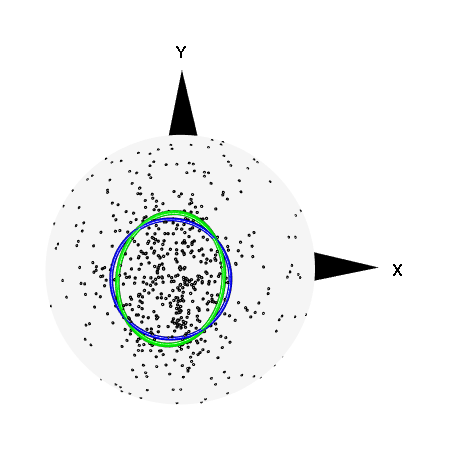}
		\caption{$\tau=0.5$}
	\end{subfigure}
	\begin{subfigure}{0.32\linewidth}
		\includegraphics[width=\linewidth]{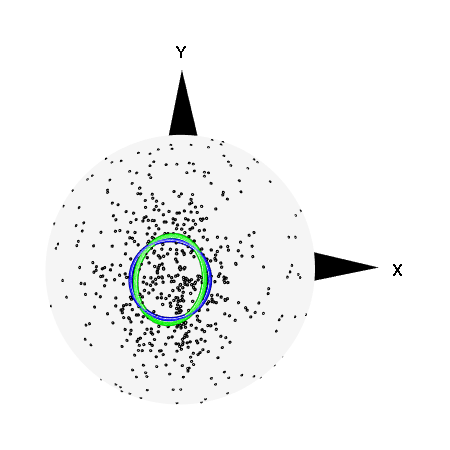}
		\caption{$\tau=0.75$}
	\end{subfigure}
	\caption{
		The empirical projection quartiles of the fibre directions $y_i$.
		The blue circle corresponds to the empirical $\tau$-depth contour $C_{\hat c_\tau \mu}$.
		The green ellipse corresponds to  $C^\cE_{\hat c^\cG_\tau\mu}$.
		The Z-direction points out of the page.
	}
	\label{fig:UHPFRC-quantiles_comparsion}
\end{figure}

As a first step, we inspect the data visually.
Figure \ref{fig:UHPFRC} shows the original and Mahalanobis-transformed fibre directions together with estimates of their densities. 
The distribution of the $y_i$ is uni-modal with empirical Fisher spherical median $\hat{\mu} =(0.029, 0.039,  0.998)^T$ which was computed by the function $mediandir()$ from the R-package $Directional$ \cite{Directional}.

Figure  \ref{fig:UHPFRC-mahala-phi} indicates that the Mahalanobis-transformed fibre directions have uniformly distributed longitudes. 
Using Watson's test applied on $\phi_{y_{i}}$ gave a p-value of less than $10^{-4}$ such that rotational symmetry is  rejected at any meaningful nominal level. 
The p-value for $\phi_{x_{i}}$ is 0.2660 such that the assumption of rotational symmetry about $\hat\mu$ is not rejected.

In Figure \ref{fig:UHPFRC-quantiles_comparsion}, we illustrate the empirical $\tau$-depth contours $C_{\hat c_\tau \mu}$ given in \eqref{eq:depth-contour} and $C^\cE_{\hat c^\cG_\tau\mu}$ given in \eqref{eq:elliptical-countour-depth} for $\tau=0.25, 0.5, 0.75$. 
The values $\minor \hat c^\cE_\tau, \major  \hat c^\cE_\tau$, and $\hat c_\tau$ are summarised in Table \ref{tab:UHPFRC-quantiles_comparsion}.
We see that $\major \hat c^\cE_\tau \lesssim \hat c^\cE_\tau  \lesssim \minor \hat c^\cE_\tau$ for all $\tau=0.25, 0.5, 0.75$. 
Thus, the shape of the underlying density seems to be slightly better fitted by an  which is elliptically symmetric distribution than by a rotationally symmetric distribution. 

\begin{table}[h]
	\centering
	\begin{tabular}{c|ccc}
		$\tau$ & 0.25 & 0.5 & 0.75 \\ 
		\hline
		$\hat{c}_\tau$          & 0.8489  & 0.9349 & 0.9729  \\
		$\minor\hat{c}^\cE_\tau$  & 0.8785  & 0.9507 & 0.9808  \\
		$\major\hat{c}^\cE_\tau$  & 0.7986  & 0.9128 & 0.9629
	\end{tabular}
	\caption{
		The empirical projection quartiles $\minor \hat c^\cE_\tau, \major  \hat c^\cE_\tau$ and $\hat c_\tau$ of the fibre directions $y_i$.
	}
	\label{tab:UHPFRC-quantiles_comparsion}
\end{table}
\subsubsection{Goodness-of-fit test and trimming}\label{sec:Goodness-of-fit test and trimming}
\begin{figure}
	\centering
	\begin{subfigure}{0.4\linewidth}
		\includegraphics[width=\linewidth]{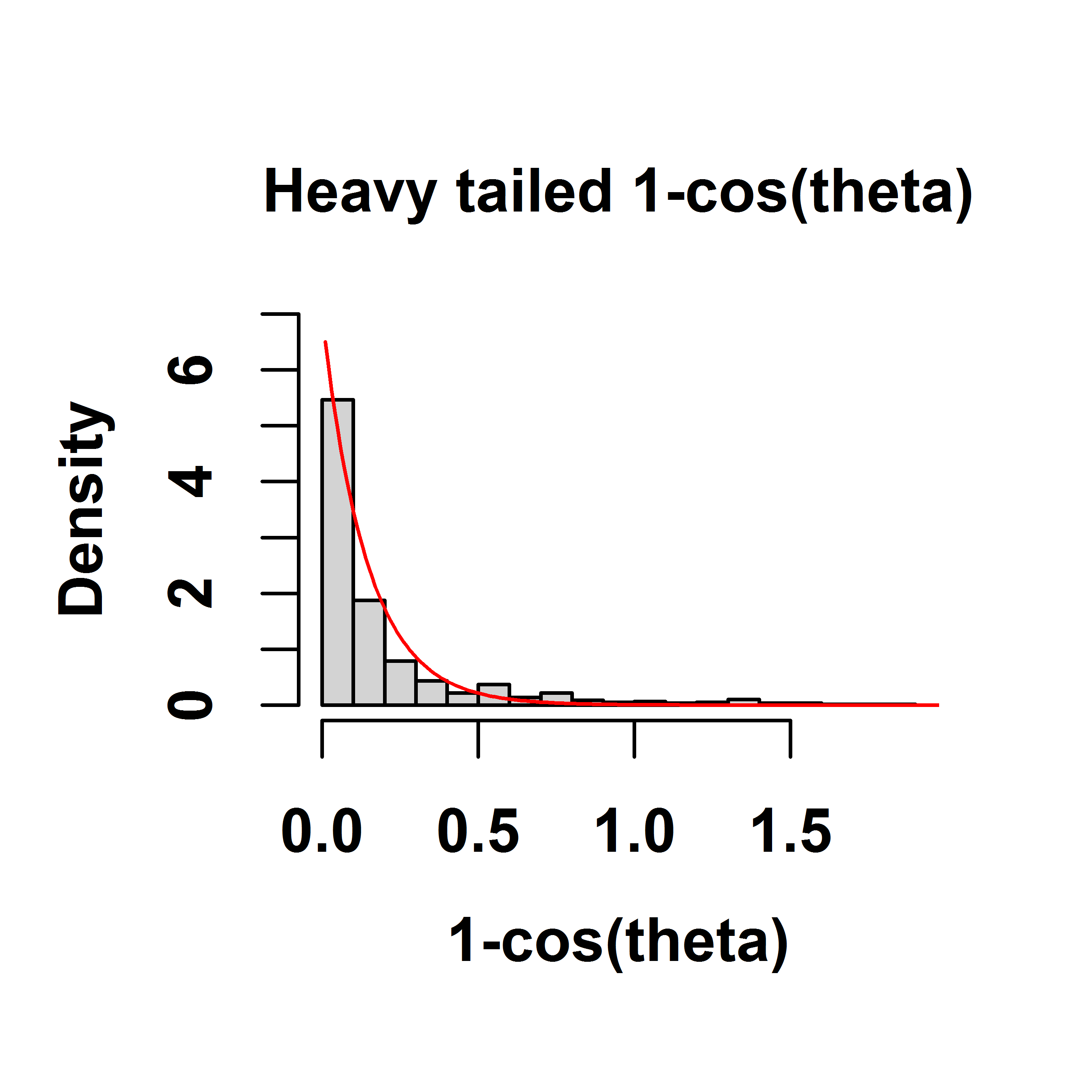}\caption{Heavy tailed $1-\cos{(\theta_i)}$ }\label{fig:UHPFRC-heavy-tail-1-costheta}
	\end{subfigure}
	\begin{subfigure}{0.4\linewidth}
		\includegraphics[width=\linewidth]{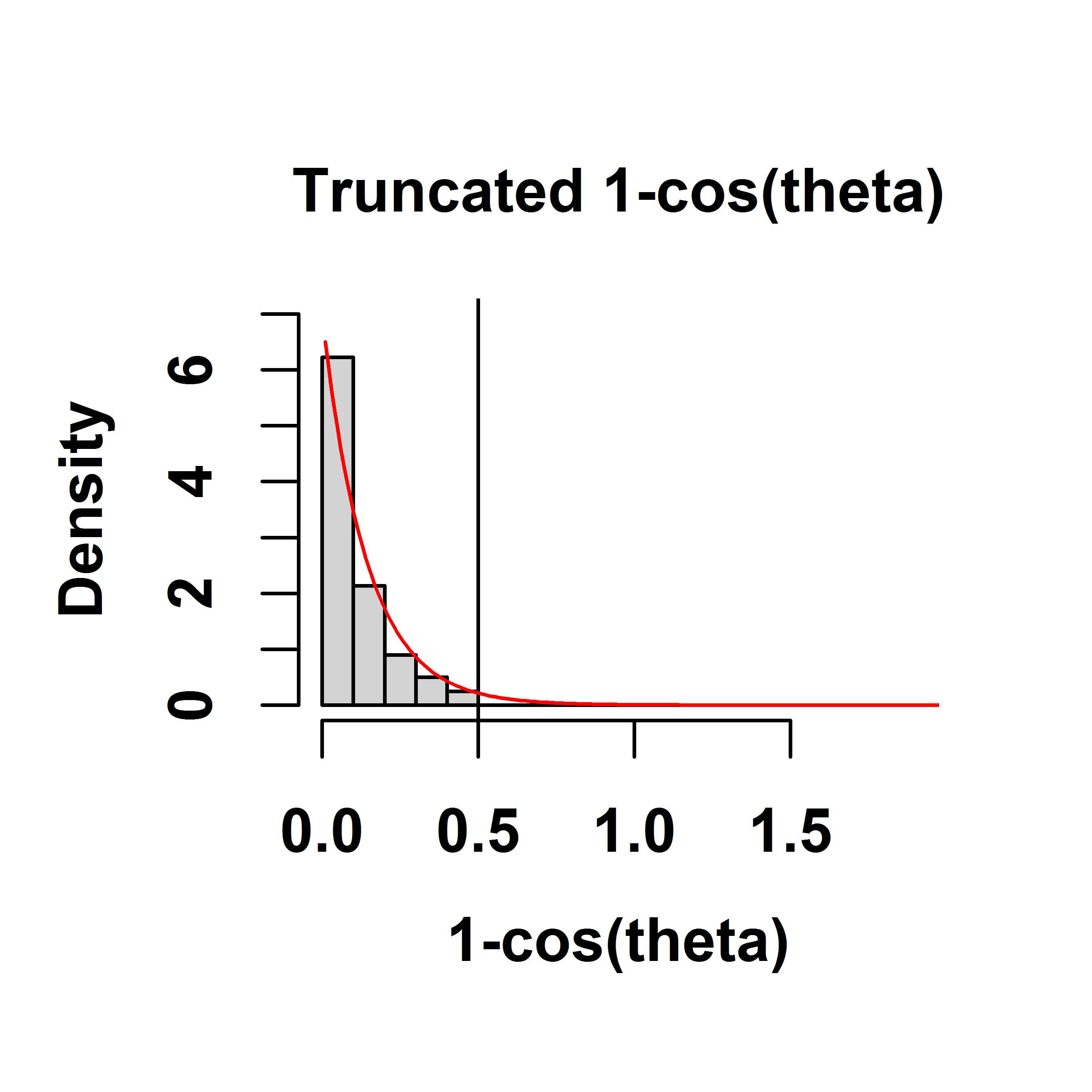}\caption{Truncated $1-\cos{(\theta_i)}$ at approx. 0.5}\label{fig:UHPFRC-trunc-1-costheta}
	\end{subfigure}
	\caption{
		Histograms of $1-\cos{(\theta_i)}$, $ i=1,\dots,n$, before (a) and after (b) truncation by trimming.
		The red line corresponds to the theoretical density of $Exp(\kappa)$, $\kappa=7$.
	}
	\label{fig:UHPFRC-heavy-tail-trunc-1-costheta}
\end{figure}
In the next step, we fit a rotationally symmetric directional distribution to the Mahalanobis transformed data $x_i$. We use a von Mises-Fisher distribution $F_0 = M_3(\mu,\kappa)$. Maximum likelihood estimation by using the function $vmf.mle()$ implemented in the R-package $Directional$ yields $\hat\kappa=6.97$. As the estimated mean direction is very close to the $Z$ axis, we consider $\mu=(0,0,1)^T$ and $\kappa=6,7,8$. We then perform the goodness-of-fit test given in \cite{Ley14} based on the projection quartiles ($\hat{c}_{0.25},\hat{c}_{0.5},\hat{c}_{0.75})$. 
The null hypothesis $H_0:F=F_0$ against $H_1:F\not=F_0$ is rejected for all three $\kappa$ values ($p< 10^{-4}$).

To analyze the reason for rejection, we convert the unit vectors $x_i$ in spherical coordinates ($\theta_i,\phi_i$).
A closer investigation of the co-latitude angle $\theta_i$ reveals that $1-\cos{(\theta_i)}$ has a heavy tail (see Figure \ref{fig:UHPFRC-heavy-tail-1-costheta}). Under the hypothesis of a von Mises-Fisher distribution, we would expect $1-\cos{(\Theta)}\sim Exp(\kappa)$ for large $\kappa$, see \cite[Eq. (4.29)]{FLE87}.
Trimming the directions $x_i$ below the $\tau$-depth contour $C_{c_\tau\mu}, \tau=0.15$, removes the heavy tail in the trimmed sample $x_i^{trim}$, see Figure \ref{fig:UHPFRC-trunc-1-costheta}. We repeat the goodness-of-fit test with $F_0=M_3(\mu,\kappa)$ with the trimmed data $x_i^{trim}$.
We chose $\mu=(0,0,1)^T$ and $\kappa=8,9,10$ because we expect a higher concentration parameter due to trimming.
The asymptotic p-values are 0.0130 ($\kappa=8$), 0.2854 ($\kappa=9$) and 0.1288 ($\kappa=10$). Thus, $\kappa=9$ provides the best fit to the trimmed  data.

In fact, the strong deviation of fibre directions from the tensile axis of almost $15 \%$ of the fibres may have favoured the cracking at this position. 

\section{Conclusion}
We extended the concept of quantiles and depth for directional data from Ley et. al. \cite{Ley14}.
Their concept provides useful geometric properties of the depth contours (such as convexity and rotational equivariance) and a Bahadur-type representation of the quantiles. 
However, a disadvantage is that rotationally symmetric depth contours are always produced, even if the underlying distribution is not rotationally symmetric \cite{Pandolfo17}. 
Our extension solves this lack of flexibility for distributions with elliptical depth contours.
The main idea was to transform the elliptical contours in the tangent space to rotationally symmetric contours, apply the results of Ley et al. \cite{Ley14}  to those, and then transform back.
In view of similarities with the classical Mahalanobis depth, our depth was called elliptical Mahalanobis depth ($EMHD_F$). 
Our results were confirmed by a Monte Carlo simulation study.
Furthermore, we introduced tools to evaluate the ellipticity of depth contours and demonstrated that our approach is the obvious choice for trimming directional data from an elliptically symmetric distribution.
We applied our quantiles and depth to analyse fibre directions in fibre-reinforced concrete.

\appendix



\bibliographystyle{unsrt} 
\bibliography{Mybibtex}

\begin{thebibliography}{10}

\bibitem{Ley14}
Ley C., Sabbah C., and Verdebout T.
\newblock {A new concept of quantiles for directional data and the angular
  Mahalanobis depth}.
\newblock {\em Electronic Journal of Statistics}, 8(1):795 -- 816, 2014.

\bibitem{Kelker70}
Kelker D.
\newblock Distribution theory of spherical distributions and a location-scale
  parameter generalization.
\newblock {\em Sankhyā: The Indian Journal of Statistics, Series A
  (1961-2002)}, 32(4):419--430, 1970.

\bibitem{Cambanis81}
Cambanis S., Huang S., and Simons G.
\newblock On the theory of elliptically contoured distributions.
\newblock {\em Journal of Multivariate Analysis}, 11(3):368--385, 1981.

\bibitem{Fang90}
Fang K., Kotz S., and Ng~K.
\newblock {\em Symmetric multivariate and related distributions}.
\newblock Chapman \& Hall, 1990.

\bibitem{Verdebout20}
Garc\'{i}a-Portugu\'{e}s E., Paindaveine D., and Verdebout T.
\newblock On optimal tests for rotational symmetry against new classes of
  hyperspherical distributions.
\newblock {\em Journal of the American Statistical Association},
  115(532):1873--1887, 2020.

\bibitem{Leong1998}
Leong P. and Carlile S.
\newblock Methods for spherical data analysis and visualization.
\newblock {\em Journal of Neuroscience Methods}, 80(2):191--200, 1998.

\bibitem{Kent82}
Kent J.
\newblock The {Fisher-Bingham} distribution on the sphere.
\newblock {\em Journal of the Royal Statistical Society. Series B
  (Methodological)}, 44(1):71--80, 1982.

\bibitem{Bahadur66}
Bahadur R.~R.
\newblock {A Note on Quantiles in Large Samples}.
\newblock {\em The Annals of Mathematical Statistics}, 37(3):577 -- 580, 1966.

\bibitem{Pandolfo17}
Pandolfo G., Paindaveine D., and Porzio G.
\newblock Distance-based depths for directional data.
\newblock {\em Canadian Journal of Statistics}, 46, 09 2017.

\bibitem{Pennec06}
Pennec X.
\newblock Intrinsic statistics on {Riemannian} manifolds: Basic tools for
  geometric measurements.
\newblock {\em Journal of Mathematical Imaging and Vision}, 25:127--154, 07
  2006.

\bibitem{Hardle03}
Härdle W. and Simar L.
\newblock {\em Applied Multivariate Statistical Analysis}.
\newblock Springer-Verlag Berlin Heidelberg, 2003.

\bibitem{Maryamh20}
Maryamh K., Hauch K., Redenbach C., and Schnell J.
\newblock Influence of production parameters on the fiber geometry and the
  mechanical behavior of ultra high performance fiber-reinforced concrete.
\newblock {\em Structural Concrete}, 22(1):361--375, 2021.

\bibitem{Mardia99}
Mardia K.~V. and P.~E. Jupp.
\newblock {\em Directional Statistics}.
\newblock Wiley, New York, 1999.

\bibitem{FLE87}
Fisher N.I., Lewis T., and Embleton B.J.J.
\newblock {\em Statistical Analysis of Spherical Data}.
\newblock Cambridge University Press, 1987.

\bibitem{Fisher85}
Fisher N.I.
\newblock Spherical medians.
\newblock {\em Journal of the Royal Statistical Society. Series B
  (Methodological)}, 47(2):342--348, 1985.

\bibitem{Tu11}
Tu~L.~W.
\newblock {\em An Introduction to Manifolds}.
\newblock Springer-Verlag, New York, 2011.

\bibitem{Fletcher08}
Fletcher P., Venkatasubramanian S., and Joshi S.
\newblock Robust statistics on {Riemannian} manifolds via the geometric median.
\newblock {\em 26th IEEE Conference on Computer Vision and Pattern Recognition,
  CVPR}, pages 1--8, 06 2008.

\bibitem{Hauberg18}
Hauberg S.
\newblock Directional statistics with the spherical normal distribution.
\newblock {\em 18th International Conference on Information Fusion}, pages
  704--711, 07 2018.

\bibitem{Angulo12}
Frontera-Pons J. and Angulo J.
\newblock Morphological operators for images valued on the sphere.
\newblock In {\em 2012 19th IEEE International Conference on Image Processing},
  pages 113--116, 09 2012.

\bibitem{Koenker05}
Koenker R.
\newblock {\em Quantile Regression}.
\newblock Cambridge University Press, 2005.

\bibitem{Kong12}
Linglong Kong and Ivan Mizera.
\newblock Quantile tomography: Using quantiles with multivariate data.
\newblock {\em Statistica Sinica}, 22(4):1589--1610, 2012.

\bibitem{Jammalamadaka01}
Rao~Jammalamadaka S. and Sengupta A.
\newblock {\em Topics in Circular Statistics}.
\newblock World Scientific, 2001.

\bibitem{Directional}
Tsagris M., Athineou G., Sajib A., Amson E., and Waldstein M.~J.
\newblock {\em Directional: A Collection of R Functions for Directional Data
  Analysis}, 2021.
\newblock R package version 5.0.

\end{thebibliography}

\end{document}